\documentclass{article}

\usepackage{graphicx, caption,subcaption}
\usepackage{amsmath, amssymb, amsthm, mathtools, dsfont}
\usepackage{color, paralist}
\usepackage[text={6.5in,9in},centering]{geometry}
\usepackage[numbers, comma, square, sort&compress]{natbib}
\usepackage{subcaption}
\usepackage{xcolor}
\usepackage{bbm}
\usepackage{enumitem}
\usepackage[hypertexnames=false,colorlinks=true,linkcolor=blue,citecolor=blue]{hyperref}
\usepackage{xurl}
\usepackage{array}

\usepackage{lineno}

\newtheorem{lemma}{Lemma}
\newtheorem{theorem}{Theorem}

\makeatletter\@addtoreset{equation}{section}\makeatother

\makeatletter

\newcommand{\tpitchfork}{\vbox{\baselineskip\z@skip
\lineskip-.52ex
\lineskiplimit\maxdimen
\m@th
\ialign{##\crcr\hidewidth\smash{$-$}\hidewidth\crcr$\pitchfork$\crcr}}}
\makeatother

\begin{document}

%
\title{Quantitative metrics for trait and identity distributions}

\author{%
 Leah Hoogstra\thanks{Division of Applied Mathematics, Brown University, Providence, RI~02912, USA}
\and 
Katherine Slyman\thanks{Department of Mathematics, Boston College, Chestnut Hill, MA~02467, USA}
\and Bj\"orn Sandstede\footnotemark[1]
}

\date{}
\maketitle


\begin{abstract}
Understanding the role of demographic diversity in group settings requires effective quantitative metrics. Intersectional feminist theory has highlighted that demographic identities can intersect in complex ways, but most metrics used to study these traits are one-dimensional. In their paper \textit{Diversity, identity, and data} (2025), Topaz et al.\ introduced two novel metrics that capture multiple aspects of demographic identities among group members: \textit{intersecting diversity} and \textit{shared identity}. We present a mathematical framework to provide probabilistic interpretations for both metrics. Using these interpretations, we prove that these two measures are anti-correlated and establish tight bounds on their possible combined values, demonstrating that there is no clear “optimal” point that maximizes both metrics. We apply these metrics in three case studies on Hollywood movies, the television show \textit{Survivor}, and a random sample of North American companies in which we explore their bounds and anti-correlation as well as their relationship to group performance in these settings. By formalizing the mathematical structure for these metrics and demonstrating their empirical relevance, we provide a foundation for researchers across the social sciences, mathematics, and related fields to more precisely quantify distributions of intersecting traits within groups and better understand their implications for group dynamics and performance.
\end{abstract}


\section{Introduction}
\label{sec:intro}
When working with a population of individuals, it is natural to consider how their traits are distributed and represented in groups.
In ecology, for instance, we might be interested in the range of different phenotypes or genotypes that are represented in microbial communities or in understanding similarities of plant communities \cite{sorenson1948method, jaccard1912distribution, wojcik_measuring_2025,zhu_trait_2017,ryabov_estimation_2022}. In the social sciences, we might be interested in the demographic profile of a specific group such as a committee or a social community and ask which identities are represented in a given group and how different identities are shared across individuals. Over the past few decades, research studies have also focused on understanding the impact of demographic diversity on group performance: while some studies have found that diverse groups can perform better on some tasks, others have found negative or no effect \cite{phillips2014diversity, mckinsey1-hunt2015diversity, mckinsey3-dolan2020diversity, herring2009does,georgeac2023business,green2024mckinsey, bermiss2024does, verwijs2023double}. The mixed results point to contextual factors mediating the relationship between diversity and performance in teams, such as the type of task, the environment, the time frame, and other characteristics of the group \cite{joshi2009role, verwijs2023double}. 

Theory for diversity's positive effect on performance is rooted in  ``cognitive resource theory", which suggests that individuals contribute diverse cognitive resources to a challenge, enabling creativity and deeper exploration when solving problems \cite{verwijs2023double, phillips2014diversity}. 
Other studies suggest that performance can increase when team members share demographic identities \cite{ryan2003stereotype, thomas1999cultural, kramer20011, milton2005identity}. 
Theory for this effect rests on the ``similarity-attraction paradigm", suggesting that shared identities may allow a team to build trust while avoiding tokenization and conflict \cite{verwijs2023double, stichman2010strength,topazdiversity}. We could then ask whether a team can optimize both of these aspects, whether there is a trade-off between creating a diverse group or focusing on groups where individuals share many demographic identities, and how we might quantify this trade-off if it is present.

Our goal is to develop, analyze, and interpret metrics that allow us to quantify the representation and distribution of individual traits in groups. Even though any such metrics will necessarily be imperfect, they may aid social scientists in their research by supplying ways to quantify trait distributions. Many such measures are already well-established \cite{budescu2012measure, roberto2015divergence, boydstun2014importance}. These broadly fall into entropic and non-entropic categories. Entropic measures are based on the information theoretic definition of entropy which logarithmically weighs probabilities to achieve useful mathematical properties, while non-entropic measures can retain direct probabilistic interpretations \cite{budescu2012measure}. In this paper, we study properties of two non-entropic quantitative indicators, namely \textit{Intersectional Diversity} ($\mathcal{D}$) and \textit{Shared Identity} ($\mathcal{S}$), that were recently introduced in \cite{topazdiversity}. Specifically, we show the following:

\begin{compactitem}
\item  Intersectional Diversity ($\mathcal{D}$) reflects the probability that two randomly chosen individuals in a given group differ in at least one trait.
\item  Shared Identity ($\mathcal{S}$) is given by the average fraction of traits shared between two randomly selected individuals.
\item The two measures $\mathcal{D}$ and $\mathcal{S}$ are anti-correlated: there is a trade-off between increasing diversity and increasing shared identities.
\item The set of possible values of $(\mathcal{D},\mathcal{S})$ is restricted to an explicit polygonal shape in the unit square.
\end{compactitem}

There are various ways to think about the possible values of $(\mathcal{D},\mathcal{S})$ depending on the goals and context of a study. A natural goal might be to get both metrics to their maximum value, but this we show is structurally impossible as they are not independent of one another. Their anti-correlation and the shape of their attainable region further highlight that there is no clear ``optimal" point that maximizes both metrics.
In this work we derive the structural constraints on attainable $(\mathcal{D},\mathcal{S})$ pairs, and the optimal values for each measure within those constraints will depend on the context in which they are applied. 

We use three real-world datasets about movie crews, the teams for the CBS reality TV show \textit{Survivor}, and financial data for a random sample of publicly traded U.S.\ companies to illustrate our results. For instance, we show how the intersectional diversity and shared identity metrics could be utilized to quantify whether programs to increase diversity have been successful (applied to CBS's pledge to create more diverse \textit{Survivor} teams), test the null hypothesis that groups are drawn randomly from a demographically representative pool (using movie crews as an example), and study the relation between the composition of leadership teams and performance (using company data).

The paper is organized as follows. Section~\ref{sec:definitions}  introduces key notation and terms, reviews the definitions of Intersecting Diversity ($\mathcal{D}$) and Shared Identity ($\mathcal{S}$), explains their probabilistic interpretations, and presents a theorem and lemma describing their bounds and anti-correlation. In Section~\ref{sec:data}, we apply these results to three case studies. We introduce the datasets, look at bounds and anti-correlation for $(\mathcal{D},{S})$, and end with the consideration of hypothesis testing. Section~\ref{sec:proofs} contains the mathematical proofs for the claims made in Section~\ref{sec:definitions}. We close with a discussion of our conclusions and paths for future work in Section~\ref{sec:discussion}. 

\section{Measures of intersectional diversity and shared identity}
\label{sec:definitions}

In this section, we introduce the terminology and notation that we use throughout this work. We discuss the two metrics introduced by Topaz et al.\ \cite{topazdiversity}, formalize their probabilistic interpretations, demonstrate their anti-correlation, and introduce their theoretical bounds. 

Since demographic identities are multifaceted, we must first establish terminology to refer to specific aspects of identity in various contexts. For the purposes of this paper, we use the term \textit{traits} or \textit{axes of identity} to refer to broad categories of identity, such as age or gender. The term \textit{values} refers to a trait's subcategories, and \textit{aggregate identities} refers to unique combinations of trait values.
For example, ``Pacific Islander man'' and ``white woman'' could be aggregate identities in a dataset with the traits ``race" and ``gender".

\begin{table}
\centering
\begin{tabular}{|m{1.75cm}|m{13cm}|}
\hline
$N$ & total number of individuals in the group \\ \hline
$T$ & total number of traits (axes of identity) being considered \\ \hline 
$v_t$ & number of values possible for trait $t$ to take on ($t \in \{ 1, 2, ... , T\})$ \\ \hline
$\boldsymbol{\mathcal{C}}$ & set of all possible aggregate identities $c=(c_1,\ldots,c_T)$ with $1\leq c_t\leq v_t$ so that $c_t$ is the value for trait $t$ for $t\in\{1,2,\ldots,T\}$\\ \hline
$C$ & total number of possible aggregate identities ($C:=|\boldsymbol{\mathcal{C}}|=v_1v_2 ... v_T)$ \\  \hline 
$p_c$ & proportion of the group that has aggregate identity $c\in \boldsymbol{\mathcal{C}}$ \\ \hline
$X$ & assigns to a pair $(i,j)$ of individuals the number $X(i,j)$ of traits they share ($X \in \{ 0, 1, ..., T \})$ \\ \hline
$P(X=x)$ & probability that two randomly selected individuals sampled with replacement share exactly $x$ traits ($x \in \{ 0, 1, ..., T \})$ \\ \hline
\end{tabular}
\caption{We summarize the mathematical notation introduced in this paper which describes various aspects of demographic identities among members of a group.}
\label{tab:notation}
\end{table}

In this work, we make the simplifying assumption that all traits have mutually exclusive values. The metrics we explore do not inherently require this, but making this assumption avoids the complexities of  
accounting for multiple choices.
Some traits naturally have this property, such as age, but in others, such as race, forcing mutual exclusivity inherently restricts the faithful representation of nuanced identities. We use this assumption purely for the initial examination of the mathematical properties of these measures. Choosing appropriate trait values is a nontrivial challenge and highly dependent on context \cite{topazdiversity}. 

To address the presence of multiple axes of identity and as a step towards accounting for the importance of intersectionality \cite{crenshaw2013demarginalizing}, Topaz et al.\ \cite{topazdiversity} introduced two new indicators, which, to our knowledge, are the first multidimensional demographic indicators intended for the social sciences. Before we define these indicators, we introduce the notation we will use. This notation is also summarized in Table~\ref{tab:notation}. Given a group of individuals and a collection of $T$ different traits, we denote by $\boldsymbol{\mathcal{C}}$ the resulting set of all possible aggregate identities determined by unique combinations of values for these $T$ traits, and by $p_c$ the proportion of group members who have aggregate identity $c\in\boldsymbol{\mathcal{C}}$. The total number of aggregate identities is denoted by $C:=|\boldsymbol{\mathcal{C}}|$. For individuals $i,j$ in our group, we let $X(i,j)$ be the number of traits shared between them: in particular, $0\leq X\leq T$.  

The first indicator introduced by \cite{topazdiversity} is the \textit{Intersecting Diversity} measure $\mathcal{D}$ defined by
\begin{align}
\label{ID_defn}
    \mathcal{D} := \frac{C}{C-1} \left( 1 - \sum_{c\in \boldsymbol{\mathcal{C}}} p_c^2 \right)
    = \frac{C}{C-1} \left( 1 - P(X = T) \right).
\end{align} 
The normalization factor $\frac{C}{C-1}$ ensures that $0\leq\mathcal{D}\leq1$, where both values can be attained; $\mathcal{D}=1$ when all traits are distributed uniformly so that $P(X = T) = \frac{1}{C}$, and $\mathcal{D}=0$ when all individuals share all traits. Since $p_c$ is the probability that a randomly selected group member has aggregate identity $c$, we see that $\mathcal{D}$ is the (normalized) probability that two randomly selected individuals belong to different aggregate identity categories when sampled with replacement. The measure $\mathcal{D}$ is an extension of the generalized variance measure  \cite{budescu2012measure}, and is its equivalent if only one axis of identity is considered \cite{topazdiversity}.

The second indicator is the \textit{Shared Identity} measure $\mathcal{S}$, which we define to be
\begin{equation}
\label{Sinf_def}
\mathcal{S} := \frac{1}{T} \sum_{t=1}^T \sum_{v=1}^{v_t} \left( \sum_{\underset{\scriptstyle c_t=v}{c\in\boldsymbol{\mathcal{C}}}} p_c\right)^2 = \frac{E(X)}{T},
\end{equation}
where $v_t$ is the number of possible values for trait $t$, and $c_t$ is the value for trait $t$ given the aggregate identity $c \in \boldsymbol{\mathcal{C}}$. The quantity $\mathcal{S}$ is equal to the expected percentage of traits shared between two randomly selected group members sampled with replacement: indeed, the term inside the brackets of \eqref{Sinf_def} is the probability that the value of trait $t$ for a randomly selected individual is $v$, and squaring this term gives the probability that $v$ is the value of trait $t$ for two randomly selected individuals sampled with replacement; summing over all $v$ gives the probability that the two sampled individuals share trait $t$, and summing the resulting expression over all $t$ gives the expected number of shared traits, as claimed. Similarly to $\mathcal{D}$, we have $0\leq\mathcal{S}\leq1$, where $\mathcal{S}=1$ is attained when all group members share all traits.

We note that \cite{topazdiversity} originally defined the Shared Identity measure
\begin{eqnarray}
\label{eqn: OG S}
\mathcal{S}_{N} &= \frac{ \text{total number of shared identity characteristics in group}}{\text{theoretical maximum number of shared identity characteristics in group}} \nonumber \\ 
&= \frac{1}{T} \frac{\sum_{i>j} X(i,j)}{\binom{N}{2}},    
\end{eqnarray}
where a \textit{shared identity characteristic} arises when two group members share the same value for a trait. As shown in the following lemma, this metric is also the expected percentage of shared traits between two randomly chosen group members, but when sampling is done \textit{without} replacement. 

\begin{lemma} 
\label{lemma:prob_int_S}
Assume a population $G$ consists of $N$ individuals, each with $T$ traits where each trait has mutually exclusive values. Let $X(i,j)$ assign to a pair of individuals $(i,j) \in G \times G$ the number $X(i,j)$ of traits they share, and let $O = \{ (i, j) \in G \times G\colon i \neq j \}$ be the set of pairs of distinct individuals in the population, so that $X|_O: O \to \{0,1,...,T\}$ is the random variable that assigns to a randomly drawn pair of distinct individuals the number of traits they share. Then
\begin{equation}\label{e:Sexp}
\mathcal{S}_N =  \frac{1}{T}  \frac{\sum_{i>j} X(i,j)}{\binom{N}{2}} = \frac{1}{T} \sum_{x=0}^T x P(X|_O=x) = \frac{E(X|_O)}{T}.
\end{equation}
\end{lemma}

Note that $X|_O$ is the same random variable as $X$ but with its domain restricted to $O$ so that it samples individuals without replacement. A proof of this lemma is provided in Section~\ref{sec:proofs}. We remark that we can analogously define a metric $\mathcal{D}_N$ based on $\mathcal{D}$ by using $P(X|_O=T)$ in \eqref{ID_defn} and adjusting the normalization constant accordingly.

The new quantity $\mathcal{S}$ aligns better with $\mathcal{D}$ since both sample with replacement. The difference between $\mathcal{S}$ and $\mathcal{S}_N$ is smaller for larger groups because $\mathcal{S} \approx \mathcal{S}_N$ when the group is large enough so that the chance of randomly selecting the same group member twice is negligible; the exact relationship between $\mathcal{S}$ and $\mathcal{S}_N$ is stated in Theorem~\ref{main_thm}.

Topaz et al.\ \cite{topazdiversity} observed that, while the values of the two indicators $\mathcal{D}$ and $\mathcal{S}_N$ can vary independently of each other, the measures appeared anti-correlated in their case studies. In Theorem~\ref{main_thm}, we prove that $\mathcal{D}$ and $\mathcal{S}$ are indeed anti-correlated quantities, describe the region where ($\mathcal{D}$, $\mathcal{S}$) pairs will occur, and provide a precise relation between $\mathcal{S}$ and $\mathcal{S}_N$. Applications of this theorem are provided in Section~\ref{sec:data}, and the proof is provided in Section~\ref{sec:proofs}. 

\begin{theorem}
\label{main_thm}
Let $\mathcal{D}$, $\mathcal{S}$, and $\mathcal{S}_N$ be defined as in \eqref{ID_defn}, \eqref{Sinf_def}, and \eqref{eqn: OG S} respectively. The indicators $\mathcal{D}$ and $\mathcal{S}$ then have the following properties:
\begin{compactenum}[(i)]
\item $\displaystyle 1-\frac{C-1}{C}\mathcal{D} \leq \mathcal{S} \leq 1-\displaystyle \frac{C-1}{TC}\mathcal{D}$;
\item $\displaystyle \mathcal{S}_\mathrm{min} :=\frac{1}{T} \sum_{t=1}^T \frac{1}{v_t} \leq \mathcal{S}$;
\item $\displaystyle \langle \nabla \mathcal{D}, \nabla \mathcal{S} \rangle \leq -\displaystyle \frac{4C}{C-1} \mathcal{S}_\mathrm{min} < 0.$
\end{compactenum}
In particular, if $T=1$, then $\mathcal{S} = 1 - \frac{C-1}{C}\mathcal{D}$. Furthermore, $\mathcal{S}=(1-\frac{1}{N})\mathcal{S}_N+\frac{1}{N}$, and each estimate of the form $a\leq\mathcal{S}\leq b$ for $\mathcal{S}$ is equivalent to the estimate $\frac{aN-1}{N-1}\leq\mathcal{S}_N\leq\frac{bN-1}{N-1}$ for $\mathcal{S}_N$.
\end{theorem}

\begin{figure}
    \centering
    \begin{subfigure}[b]{0.45\textwidth}
        \centering
        \includegraphics[scale=1]{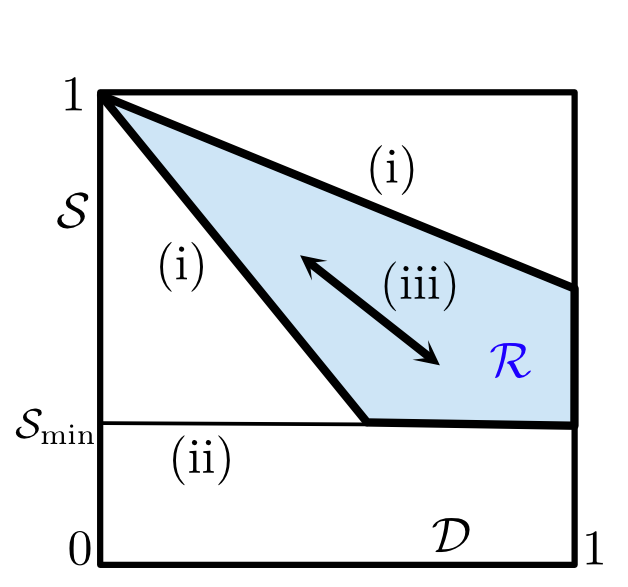} 
        \caption{}
        \label{subfig:theorem1viz}
    \end{subfigure}
    \hspace{5mm}
    \begin{subfigure}[b]{0.45\textwidth}
        \centering
        \includegraphics[scale=1]{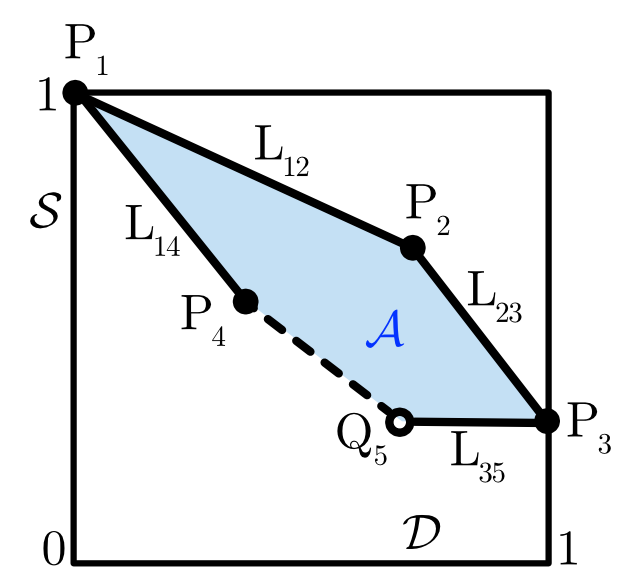} 
        \caption{}
        \label{subfig:lemma2viz}
    \end{subfigure}
    \caption{Panel~(a) illustrates that each attainable pair $(\mathcal{D},\mathcal{S})$ must lie in the polygonal region $\mathcal{R}$, where the bounds (i) and (ii) are given in Theorem~\ref{main_thm}, while (iii) indicates that the data are anti-correlated inside this region. For the case of precisely two traits ($T=2$), panel~(b) illustrates the boundaries of the attainable region $\mathcal{A}$ established in Lemma~\ref{lemma:admissible}: each point on the solid line segments can be realized by a pair $(\mathcal{D},\mathcal{S})$; the dotted line is a conjecture.}
    \label{fig:theorem-visuals}
\end{figure}

Parts $(i)$ and $(ii)$ of Theorem~\ref{main_thm} imply that any pair $(\mathcal{D}, \mathcal{S})$ must lie in the polygonal region
\[
\mathcal{R} := \left\{ (x_1,x_2)\in[0,1]^2\colon
\frac{C-1}{TC} x_1 \leq 1-x_2 \leq \frac{C-1}{C} x_1,\;
x_2\geq\frac{1}{T} \sum_{t=1}^T \frac{1}{v_t}\right\}
\]
for a given choice of traits and trait values, and part $(iii)$ shows that these quantities are anti-correlated and provides an upper limit for the correlation coefficient. The bounds depend only on the values of $v_t$, $C$, and $T$, all of which are known in advance; Figure~\ref{subfig:theorem1viz} shows an example configuration. The resulting bounds are not externally imposed, but rather arise naturally from the definitions of the metrics.
An intuitive argument for the anti-correlation of $\mathcal{D}$ and $\mathcal{S}$ follows from their interpretations. Given two randomly sampled individuals from a group, $\mathcal{D}$ captures the likelihood that they have different aggregate identities, whereas $\mathcal{S}$ captures how many traits they are likely to share. Though these interpretations are not exact inversions of each other, they have opposite effects. 

The bounds given in Theorem~\ref{main_thm} are simple to calculate and apply to any number of traits $T$, but numerical work shown in Section~\ref{sec:bounds_results} indicates that these bounds are not tight. Our final result provides an explicit expression for the set of all attainable pairs $(\mathcal{D},\mathcal{S})$ for the case of exactly two traits ($T=2$), and a visual example is shown in Figure~\ref{subfig:lemma2viz}.

\begin{lemma} 
\label{lemma:admissible}
Consider two traits with mutually exclusive values where $v_1\leq v_2$. The set
\[
\mathcal{A}:=\{(\mathcal{D},\mathcal{S})\colon \exists \mbox{ probability measure }p:\mathcal{C}\to[0,1] \mbox{ that yields } (\mathcal{D},\mathcal{S})\}
\]
has the form given in Figure~\ref{subfig:lemma2viz} where
\begin{align*}
P_1 = \left(0,1\right), \;
P_2 = \left(\frac{v_1(v_2-1)}{C-1}, \frac{v_2+1}{2v_2} \right), \;
P_3 = \left(1, \frac{v_1+v_2}{2C} \right), \; \\
P_4 = \left(\frac{(v_1-1)v_2}{C-1} , \frac{1}{v_1} \right), \;
Q_5 = \left(q_5, \frac{v_1+v_2}{2C} \right)
\end{align*}
for a unique $q_5$ with $\frac{(v_1-1)v_2}{C-1}<q_5<1$ and the curves $L_{ij}$ are line segments connecting $P_i$ to $P_j$. The lines $\mathcal{S}=1-\frac{C-1}{2C}\mathcal{D}$ and $\mathcal{S}=1-\frac{C-1}{C}\mathcal{D}$, which describe the upper and lower bounds for $\mathcal{S}$ in Theorem~\ref{main_thm}(i), contain the line segments $L_{12}$ and $L_{14}$, respectively.
\end{lemma}

The region $\mathcal{A}$ consists of every point which is theoretically attainable; that is, we can construct a group with proportions $p_c$ to attain any $(\mathcal{D}, \mathcal{S})$ value inside or on the boundary of $\mathcal{A}$. Conversely, $(\mathcal{D}, \mathcal{S})$ values outside $\mathcal{A}$ are unattainable for any group.


\section{Case studies}
\label{sec:data}
We use three real-world datasets to illustrate the bounds for $(\mathcal{D},\mathcal{S})$ and ($\mathcal{D}, \mathcal{S}_N$) pairs as well as the inherent trade-offs between $\mathcal{D}$ and $\mathcal{S}$. We provide case studies to outline how $\mathcal{D}$ and $\mathcal{S}$ can be used to answer questions about individual datasets. Specifically, we use $(\mathcal{D},\mathcal{S})$ pairs to  evaluate a policy change to increase diversity, test whether teams are drawn randomly from an underlying pool, and to check for correlation between the composition of teams and performance indicators. A replication package is available at \href{https://github.com/sandstede-lab/quantitative-trait-distributions}{https://github.com/sandstede-lab/quantitative-trait-distributions}. 

\subsection{Datasets}\label{sec:datasets}

\paragraph{Movie dataset}

The movie dataset, also analyzed in \cite{topazdiversity}, consists of teams of 10 key contributors for each of the 100 highest grossing U.S.\ films of 2018 and 2019 as identified by the Internet Movie Database (IMDB). This list contains presumptive gender, race/ethnicity, and job attribute for the majority of the contributors, and we remove movies listing fewer than 10 contributors, those with missing demographic data for any of the key contributors, and duplicate entries that appear in both 2018 and 2019.
After this preprocessing step, the final list contains 180 films with 1800 principal roles and 1496 unique individuals. The identity traits available are binary gender (with the two trait values female and male) and race/ethnicity (with the six values African American, Asian, Caucasian, Latinx, Middle Eastern, and Pacific Islander). We use trait values as they were assigned in \cite{topazdiversity}.

\paragraph{Survivor dataset}

The CBS television show \textit{Survivor} is a reality-competition program in which 16 to 20 contestants are sent to a remote location, compete in challenges, and vote to eliminate other contestants in order to win a cash prize.
A dataset that includes details on the cast demographics, challenge outcomes, and voting history for each season is available through the CRAN survivoR package \cite{survivor_cran}. During the team portion of the game, there are two to four ``tribes" of equal size to begin, but as players are voted out after losing immunity challenges, tribe sizes can vary. We filter to only include situations where there are two or more tribes consisting of four or more teammates each. Our analysis encompasses the 46 completed seasons between Spring 2000 and Spring 2024, which includes 697 unique contestants. The identity traits we use are age (with the five values 24 and younger, 25-34, 35-44, 45-54, and 55 and older), gender (with the three values female, male, and non-binary), and race/ethnicity (with the six values African American, Asian, Caucasian, Latinx, Native American, and multiracial). The gender and race/ethnicity traits were encoded in \cite{survivor_cran}, except for multiracial, which we created for the nine contestants who had more than one race/ethnicity trait listed. Age bins were chosen to be approximately equal in size and represent approximate generational differences, and varying the cutoffs slightly produced similar results. 

\paragraph{North American company dataset}

McKinsey's \textit{Diversity Matters} studies \cite{mckinsey1-hunt2015diversity, mckinsey2-hunt2018delivering, mckinsey3-dolan2020diversity} are often cited as support for the case that publicly traded U.S.\ companies perform on average better when they have more diverse executive leadership teams \cite{green2024mckinsey}. Green and Hand \cite{green2024mckinsey} aimed at a quasi-replication of these results using their own dataset of publicly traded North American companies and found that they could not replicate the findings from the McKinsey studies. We focus here on an anonymized random sample of the more extensive dataset used in \cite{green2024mckinsey, bermiss2024does} provided to us by the authors. From the random sample, we remove companies with incomplete race/ethnicity and gender data for anyone in their listed leadership teams and those with four or fewer people on their leadership teams, leaving us with 292 companies. We also remove those that went bankrupt or were acquired since 2019, leaving us with 213 companies to investigate. We exclude bankrupt and acquired companies as it is unclear how to appropriately assign their EBIT margin values. We pulled annual financial data from Compustat on July 22, 2024 for the fiscal years 2019-2023. As in  \cite{green2024mckinsey}, we used industry-adjusted EBIT margin as the performance metric, aligning with the Fama-French 12 industry categorizations. After successfully replicating the metrics in \cite{green2024mckinsey} for the Fortune~500 companies in 2019, we averaged each company's EBIT margins over the four years and adjusted them by subtracting the average of the industry median for all North American companies over those four years. The identity traits we use, both presumptive, are binary gender (the two values female and male) and race/ethnicity (the five values American Indian or Alaska Native; Asian, Native Hawaiian, or Pacific Islander; Black; Hispanic or Latino; and white). We use trait values as they were assigned in \cite{green2024mckinsey}.

\subsection{Bounds and anti-correlation}\label{sec:bounds_results}

In this section, we use the movie and \textit{Survivor} datasets to demonstrate the applications of Theorem~\ref{main_thm} and Lemma~\ref{lemma:admissible}. The data adhere to the polygonal bounds established in Theorem~\ref{main_thm}(i)-(ii) and Lemma~\ref{lemma:admissible} and are anti-correlated as indicated in Theorem~\ref{main_thm}(iii).

\begin{figure}
    \centering
    \includegraphics[width = 0.6\textwidth]{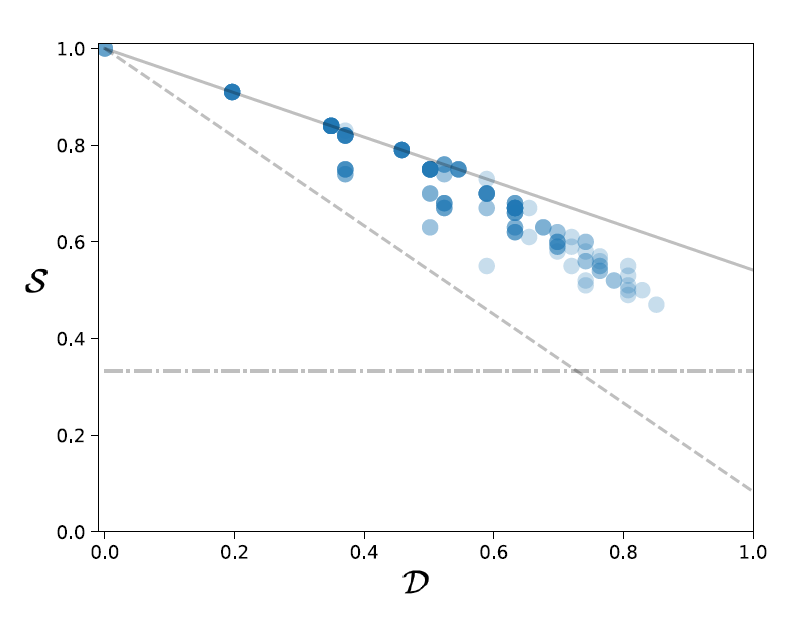}
    \caption{The circles represent the $(\mathcal{D},\mathcal{S})$ pairs of the movie dataset with $T=2$ and $v_1 = 2, v_2 = 6$ ($C=v_1v_2=12$). The circles are partially transparent, so darker circles represent more movie teams. The solid line is the upper bound $\mathcal{S}=1 - \frac{C-1}{TC} \mathcal{D} = 1 - \frac{11}{24}\mathcal{D}$, the dashed line is the lower bound $\mathcal{S}=1 - \frac{C-1}{C}\mathcal{D} = 1 - \frac{11}{12}\mathcal{D}$, and the dot-dashed line represents $\mathcal{S}_{\text{min}}= \frac{1}{3}$. }
    \label{fig:M_bounds}
\end{figure}

The movie data has $T=2$, $v_1 = 2$, and $v_2 = 6$ so that $C = 12$. Since each team has 10 members, we have a fixed relationship between $\mathcal{S}$ and $\mathcal{S}_N$: $\mathcal{S} = (1 - \frac{1}{10})\mathcal{S}_N + \frac{1}{10}$.  Figure~\ref{fig:M_bounds} shows the resulting bounds from Theorem~\ref{main_thm} together with the 180 data points, one for each team. The data are contained inside the admissible polygon $\mathcal{R}$ from Theorem~\ref{main_thm} and are anti-correlated. 

\begin{figure}
    \centering
    \begin{subfigure}[b]{0.45\textwidth}
        \centering
        \includegraphics[width=\textwidth]{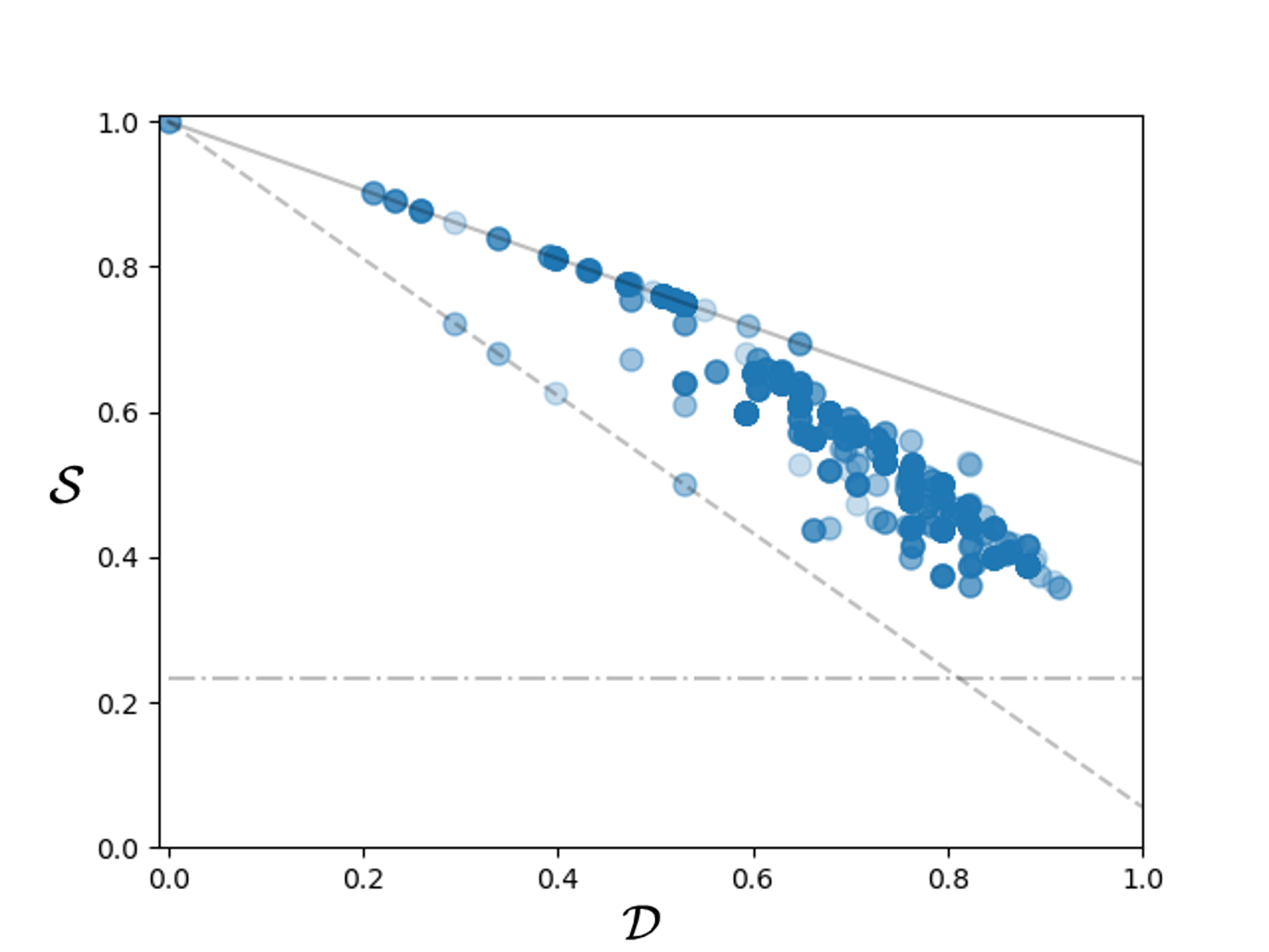} 
        \caption{}
        \label{subfig:Surv_S_bounds}
    \end{subfigure}
    \hspace{5mm}
    \begin{subfigure}[b]{0.45\textwidth}
        \centering
        \includegraphics[width=\textwidth]{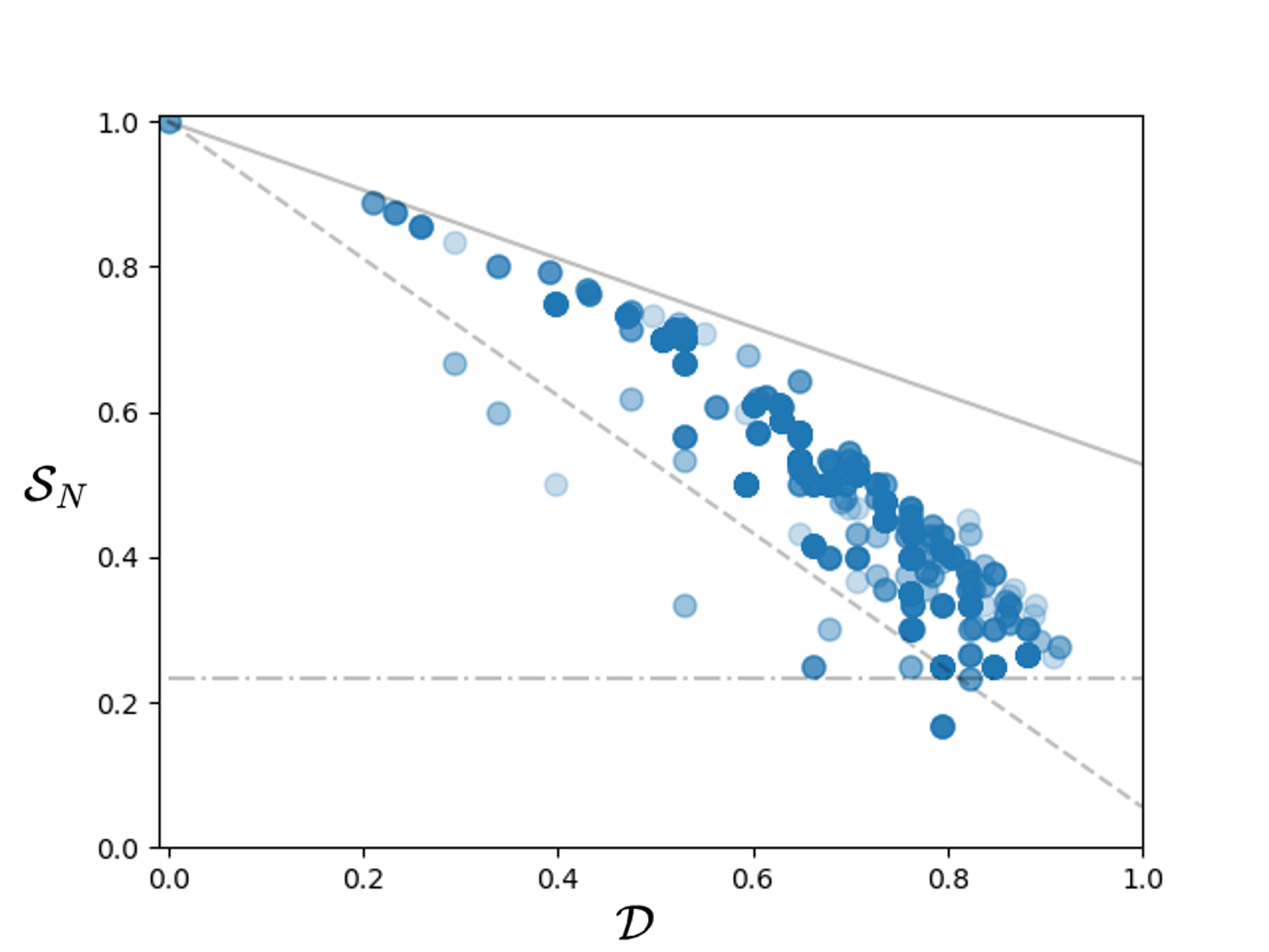} 
        \caption{}
        \label{subfig:Surv_SN_bounds}
    \end{subfigure}
    \caption{Shown are the data of the \textit{Survivor} dataset for $T=2$ and $v_1 = 3, v_2 = 6$ ($C = 18$) with $(\mathcal{D},\mathcal{S})$ in panel~(a) and  $(\mathcal{D},\mathcal{S}_N)$ in panel~(b). The circles are partially transparent, and darker circles represent more teams.
    Both panels display the bounds defined for $(\mathcal{D}, \mathcal{S})$ pairs. The solid line is the upper bound $\mathcal{S}=1 - \frac{C-1}{TC} \mathcal{D} = 1 - \frac{17}{36}\mathcal{D}$, 
    the dashed line is the lower bound 
    $\mathcal{S}=1 - \frac{C-1}{C}\mathcal{D} = 1 - \frac{17}{18}\mathcal{D}$,  
    and the dot-dashed line represents 
    $\mathcal{S}_{\text{min}} = \frac{1}{T} \sum_{t = 1}^T \frac{1}{v_t} = \frac{1}{4}$. 
    For Panel~(b), the $(\mathcal{D}, \mathcal{S})$ bounds capture the $(\mathcal{D}, \mathcal{S}_N)$ data for all tribes of size seven and larger.}
    \label{fig:S_bounds}
\end{figure}

For the \textit{Survivor} data, the results for using two or three traits are similar; here, we present the case $T=2$ with the traits gender and race/ethnicity for which $v_1 = 3$, $v_2 = 6$, and $C = 18$. Figure~\ref{subfig:Surv_S_bounds} contains the points $(\mathcal{D}, \mathcal{S})$ associated with each tribe in the dataset, and we see that these data points lie inside the admissible polygonal region $\mathcal{R}$ and are anti-correlated in line with the results in Theorem~\ref{main_thm}. Next, in Figure~\ref{subfig:Surv_SN_bounds} we consider the data points $(\mathcal{D}, \mathcal{S}_N)$ for each tribe in the dataset. The bounds for $(\mathcal{D}, \mathcal{S}_N)$ provided in Theorem~\ref{main_thm} depend on $N$, which varies by tribe, and so for small values of $N$, the bounds for $(\mathcal{D}, \mathcal{S}_N)$ and $(\mathcal{D}, \mathcal{S})$ may differ significantly. This is visible in Figure~\ref{subfig:Surv_SN_bounds} where we plot $(\mathcal{D}, \mathcal{S}_N)$ across tribes together with the admissible polygon $\mathcal{R}$ for $(\mathcal{D}, \mathcal{S})$. We see that the data points $(\mathcal{D}, \mathcal{S}_N)$ do not all lie in $\mathcal{R}$, and the tribes outside the boundaries all have six or fewer team members. This underscores the inherent dependence of $\mathcal{S}_N$ on the individual tribe size $N$, resulting in a larger deviation from $\mathcal{S}$ for smaller tribes. For example, duplicating the group while maintaining the same proportions $p_c$ can inflate $\mathcal{S}_N$, but $\mathcal{S}$, like $\mathcal{D}$, only depends on proportional representation within the group. This exploration highlights that while $\mathcal{S}_N$ can be a useful tool when $N$ is fixed across a dataset, $\mathcal{S}$ allows for a direct comparison across differently-sized groups. For the remainder of our analyses, we focus on $(\mathcal{D}, \mathcal{S})$.

\begin{figure}
    \centering
     \begin{subfigure}[b]{0.49\textwidth}
        \centering
        \includegraphics[width=\textwidth]{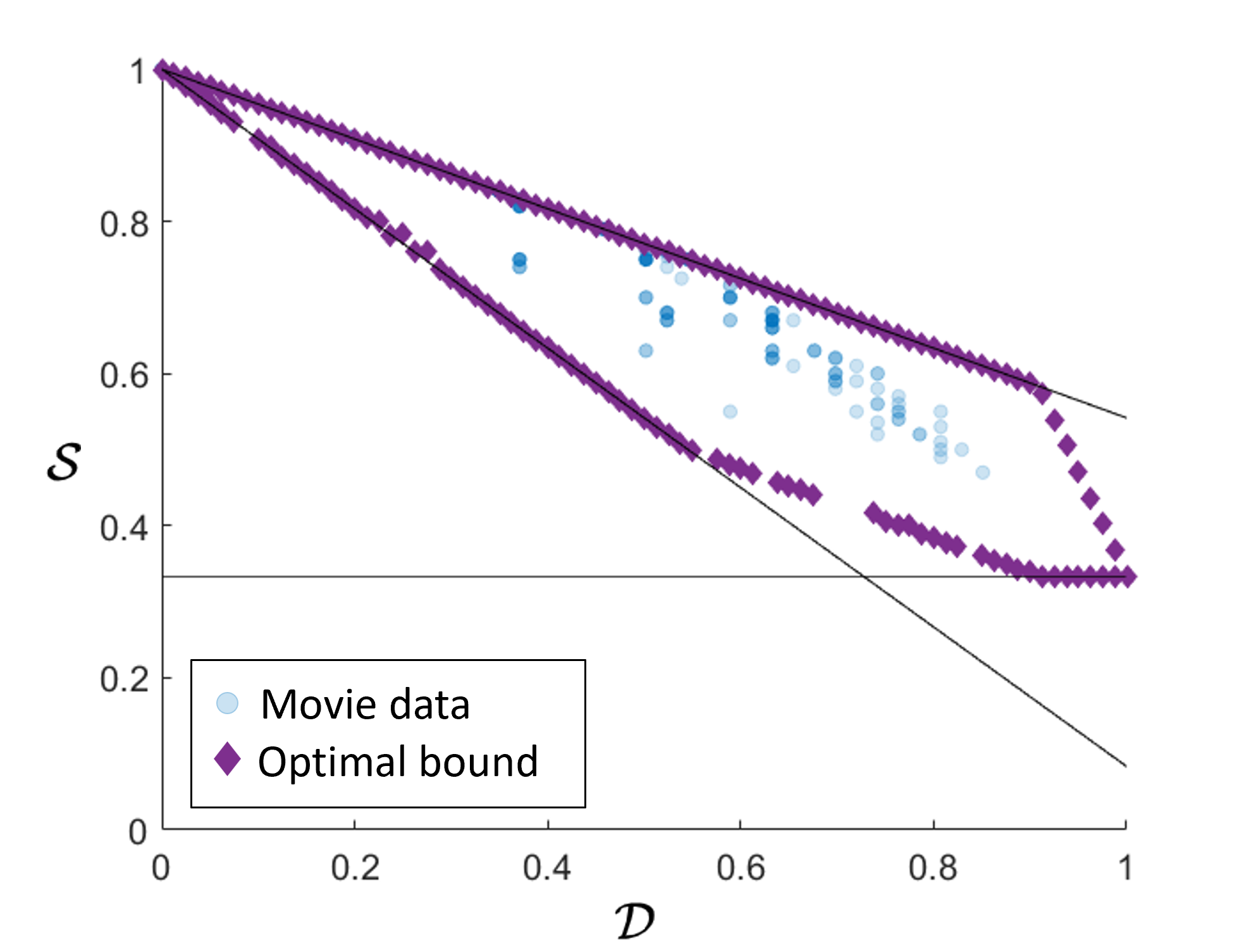} 
        \caption{}
        \label{subfig:movie_numerics}
    \end{subfigure}
    \begin{subfigure}[b]{0.49\textwidth}
        \centering
        \includegraphics[width=\textwidth]{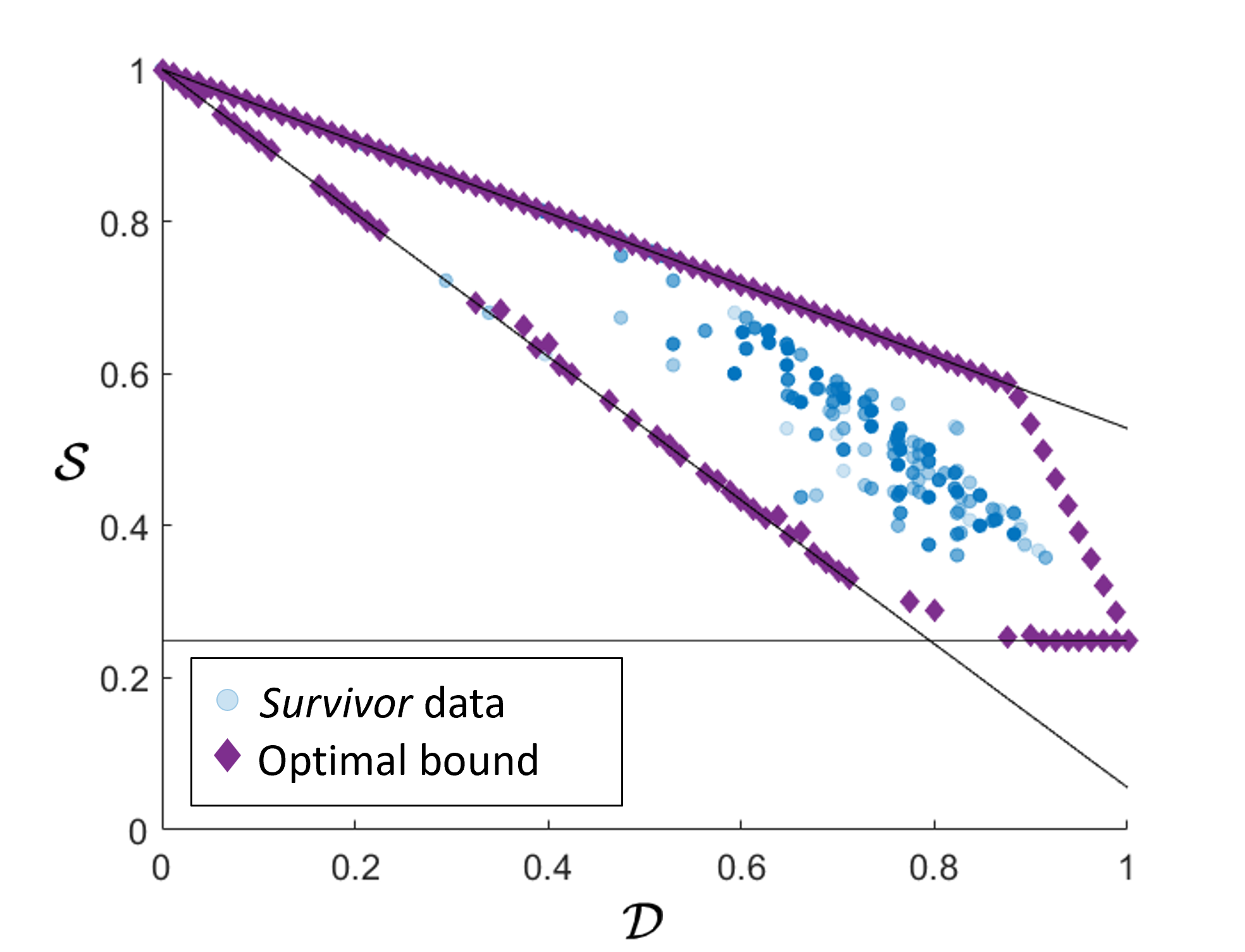} 
        \caption{}
        \label{subfig:survivor_numerics}
    \end{subfigure}
    \caption{The panels show the optimal bounds for $(\mathcal{D}, \mathcal{S})$ obtained by numerical optimization (purple diamonds) together with the bounds proved in Theorem \ref{main_thm} (gray lines) for two traits with $v_1 = 2$ and $v_2 = 6$ in panel (a) and  $v_1 = 3$ and $v_2 = 6$ in panel~(b). Panel~(a) also contains the movie data and panel~(b) the \emph{Survivor} data (both in blue circles) for comparison. All data lie inside or on the polygon outlined by the optimized boundary points.}
    \label{fig:numerical_bounds}
\end{figure}

Finally, we examined the boundaries of the attainable region $\mathcal{A}$ for two different cases with two traits ($T=2$), first with $v_1 = 2, v_2 = 6$ (corresponding to the movie dataset) and second with $v_1 = 3, v_2 = 6$ (corresponding to the \textit{Survivor} dataset). The polygonal region $\mathcal{A}$ given in Lemma~\ref{lemma:admissible} provides the optimal bounds for $(\mathcal{D}, \mathcal{S})$ pairs in these two cases. In Figure~\ref{fig:numerical_bounds}, we plot the bounds derived in Theorem~\ref{main_thm} overlaid with the movie and \textit{Survivor} data for their respective cases. We also show the results of the numerical optimization (separately maximizing and minimizing) of $\mathcal{S}$ given $\mathcal{D}$ for $0\leq\mathcal{D}\leq1$, using the fact that both metrics are quadratic forms and have the constraint $\sum_{c \in \mathcal{C}} p_c = 1$. The optimal numerical bounds agree with the boundary of the attainable region $\mathcal{A}$ described in Lemma~\ref{lemma:admissible}, and they also either match or improve upon the bounds provided in Theorem~\ref{main_thm}.

\subsection{Hypothesis testing}

\paragraph{Quantifying changes in team diversity}
\label{sec:changes}

The \textit{Survivor} dataset allows for the richest exploration of team composition over differing time intervals. In November 2020, multiple sources \cite{Deadline, Variety} released articles that CBS set a diversity target for their reality casts to be at least 50\% BIPOC (Black, Indigenous, (and) People of Color), which was to be implemented for the 2021-2022 broadcast season. To study this target, we calculate $(\mathcal{D},\mathcal{S})$ using the three traits gender, race/ethnicity, and age for the full cast for each season of \textit{Survivor}. We treat a season's full cast as a single group in this analysis, regardless of which tribes the cast members joined.

\begin{figure}
  \centering
  \begin{subfigure}[c]{0.6\textwidth}
    \centering
    \includegraphics[width=\linewidth]{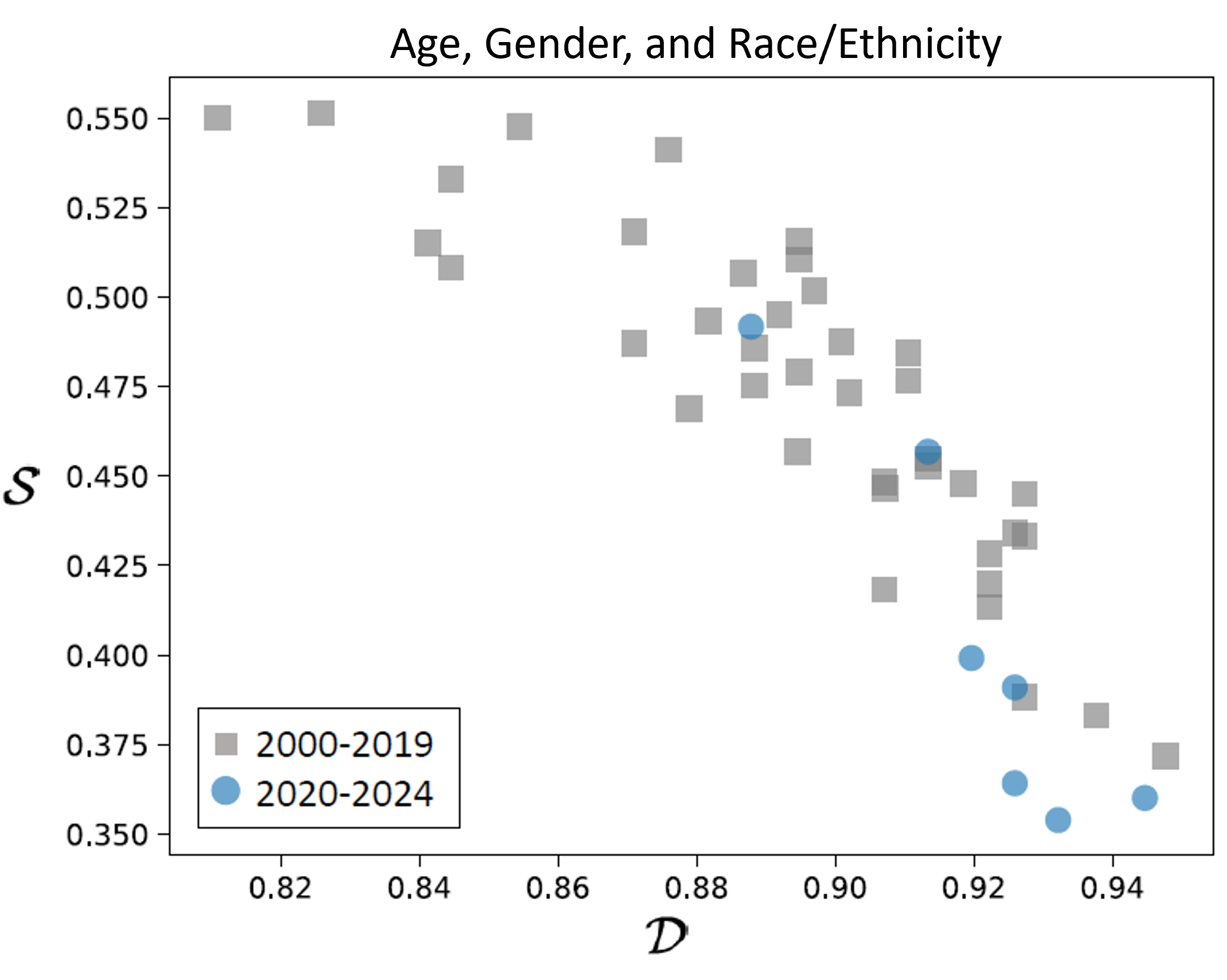}  
    \caption{} 
    \label{fig:S_3D}
  \end{subfigure}
  \begin{subfigure}[c]{0.3\textwidth}
    \centering
    \includegraphics[width=\linewidth]{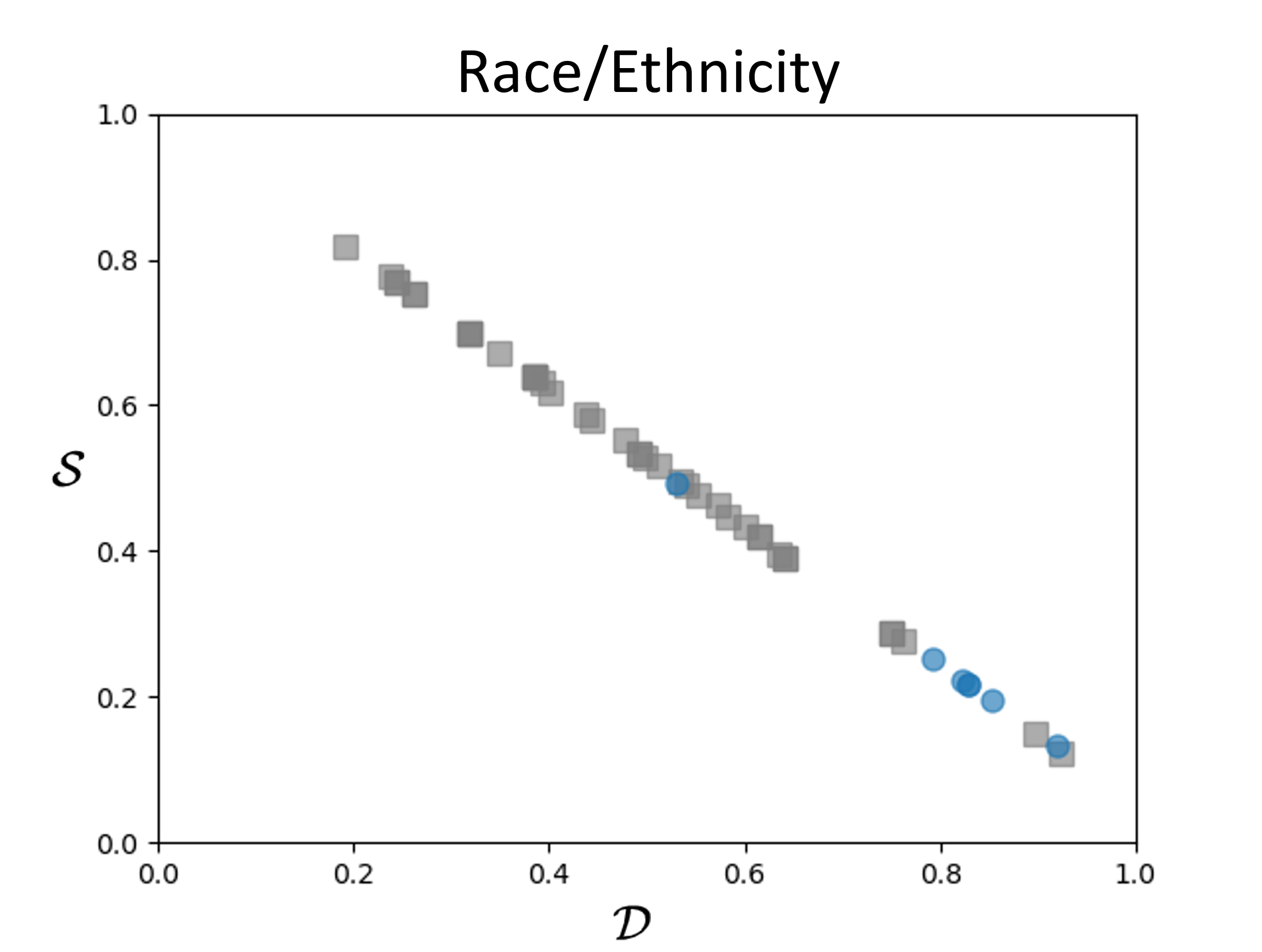}     
    \caption{}
    \label{fig:1D_race}
    \includegraphics[width=\linewidth]{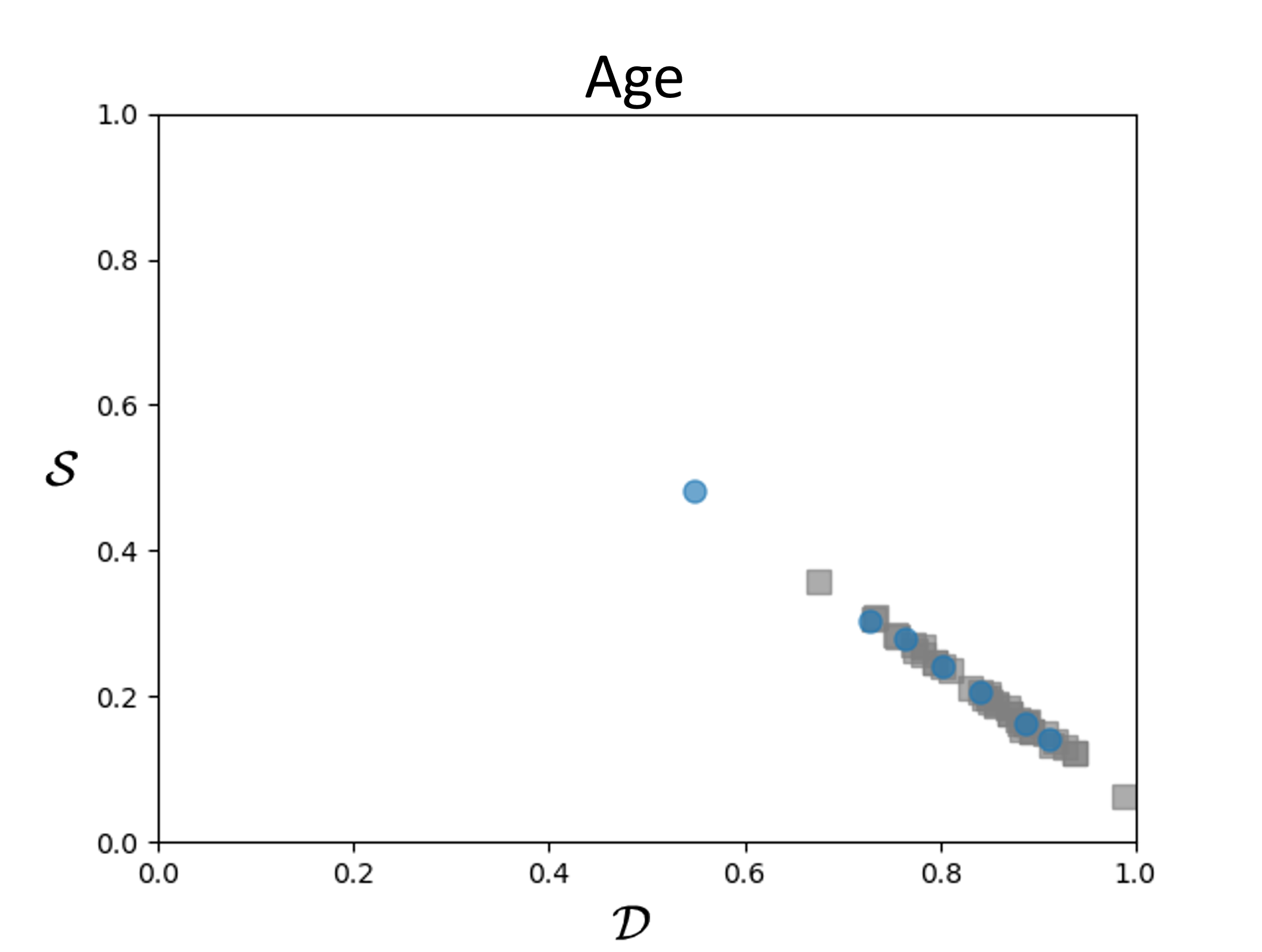}  
    \caption{}
    \label{fig:1D_age}
  \end{subfigure}
  \caption{(a) $(\mathcal{D},\mathcal{S})$ pairs for the traits age, gender, and race/ethnicity ($T=3$). Blue circles represent the seasons of \textit{Survivor} from 2020 to 2024 while gray squares denote the prior seasons. All markers are partially transparent, so darker markers indicate more seasons. Note this figure has been zoomed in on and not shown in the standard $[0,1] \times [0,1]$ box. (b) $(\mathcal{D},\mathcal{S})$ pairs for $T=1$ considering race/ethnicity only. There is an increase in $\mathcal{D}$ in years 2020-2024 (blue circles, partially transparent) compared to most prior years (gray squares, partially transparent). (c) $(\mathcal{D},\mathcal{S})$ pairs for $T=1$ considering age only. There is a decrease in $\mathcal{D}$ in years 2020-2024 (blue circles, partially transparent) compared to most prior years (gray squares, partially transparent). In panels (b)-(c), the values lie on a line since $\mathcal{D}=1-\mathcal{S}$ when $T=1$.}
  \label{fig:CBS_composite}
\end{figure}

In Figure~\ref{fig:S_3D}, we see that the cast of seasons beginning in 2020 is more diverse than those of prior seasons. Further insight can be gained by analyzing each trait individually. In Figures~\ref{fig:1D_race} and~\ref{fig:1D_age}, we plot $(\mathcal{D},\mathcal{S})$ for the single traits of race/ethnicity and age, respectively, and note that these data lie on a line since $\mathcal{S} = 1 - \frac{C-1}{C}\mathcal{D}$ when $T=1$. We see  race/ethnicity has a wider spread than age, and thus may influence the three-dimensional $\mathcal{D}$ metric more than age does. We conclude that, though the composite $\mathcal{D}$ score is notably improved after 2020 due to increased racial diversity, age diversity among the cast members actually decreased. In particular, traits can affect the overall scores very differently.

In the \emph{Survivor} dataset and others, we generally observed that the spread of ($\mathcal{D}, \mathcal{S})$ pairs is influenced by the variation of each trait both between groups and across all individuals in the dataset. To illustrate how the variation across all individuals contributes, consider a $\mathcal{D}$ score calculated from a single trait ($T=1$, as in Figures~\ref{fig:1D_age} or \ref{fig:1D_race}) and then re-calculated by adding a binary gender trait ($T=2$). If the gender trait has an even 50/50 split between its two values across the dataset, it contributes to high $\mathcal{D}$ scores, but $\mathcal{D}$ and $\mathcal{S}$ largely retain the linear relationship they had for $T=1$. However, if there is an 85/15 ratio between two gender values across the dataset, we observe a wider and less linear spread across teams. A full analysis of each trait's influence on the $\mathcal{D}$ and $\mathcal{S}$ scores is an important area for further research.

\paragraph{Representative group composition}

\label{sec:groupcomp}

A common use of group composition measures is to assess whether groups are representative of an underlying population. We can assess whether a set of teams are representative of a given population by asking if the distribution of their $(\mathcal{D}, \mathcal{S})$ scores was likely generated by composing teams randomly from the underlying population. 
We use the movie data to illustrate this type of analysis where the underlying population consists of all individuals appearing in the movie dataset. We sought to determine whether the distribution of $\mathcal{D}$ and $\mathcal{S}$ across the 180 teams could arise through random sampling of the population of available crew members, that is, through random assignments of individuals to movies. Movie teams are carefully crafted, so 
our hypothesis is that the observed ensemble of 180 movie crews is not the result of a random draw. 

\begin{figure}
\centering
\begin{subfigure}[b]{0.45\textwidth}
\includegraphics[width=\textwidth]{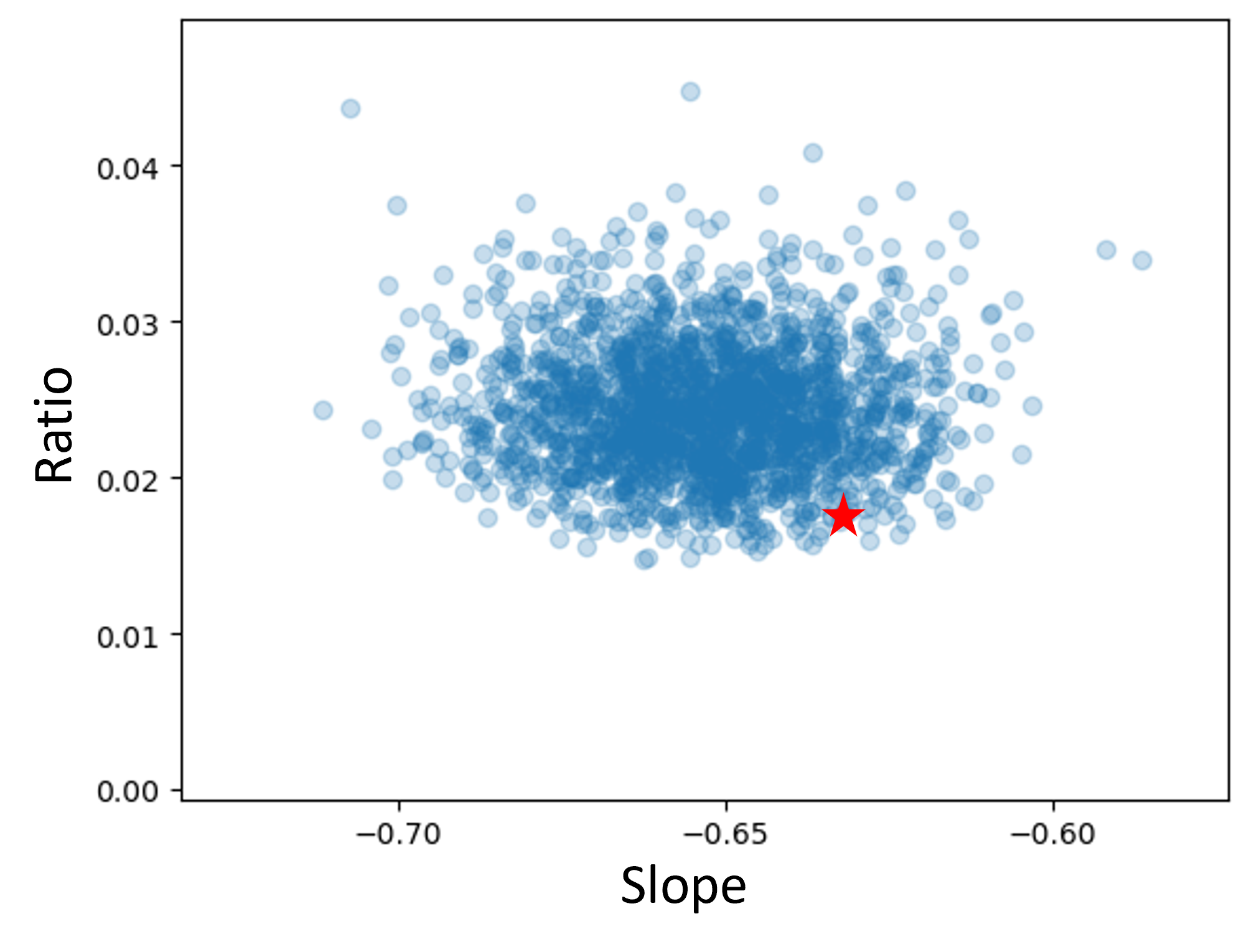}
\caption{}
\end{subfigure}
\hspace{5mm}
\begin{subfigure}[b]{0.45\textwidth}
\includegraphics[width=\textwidth]{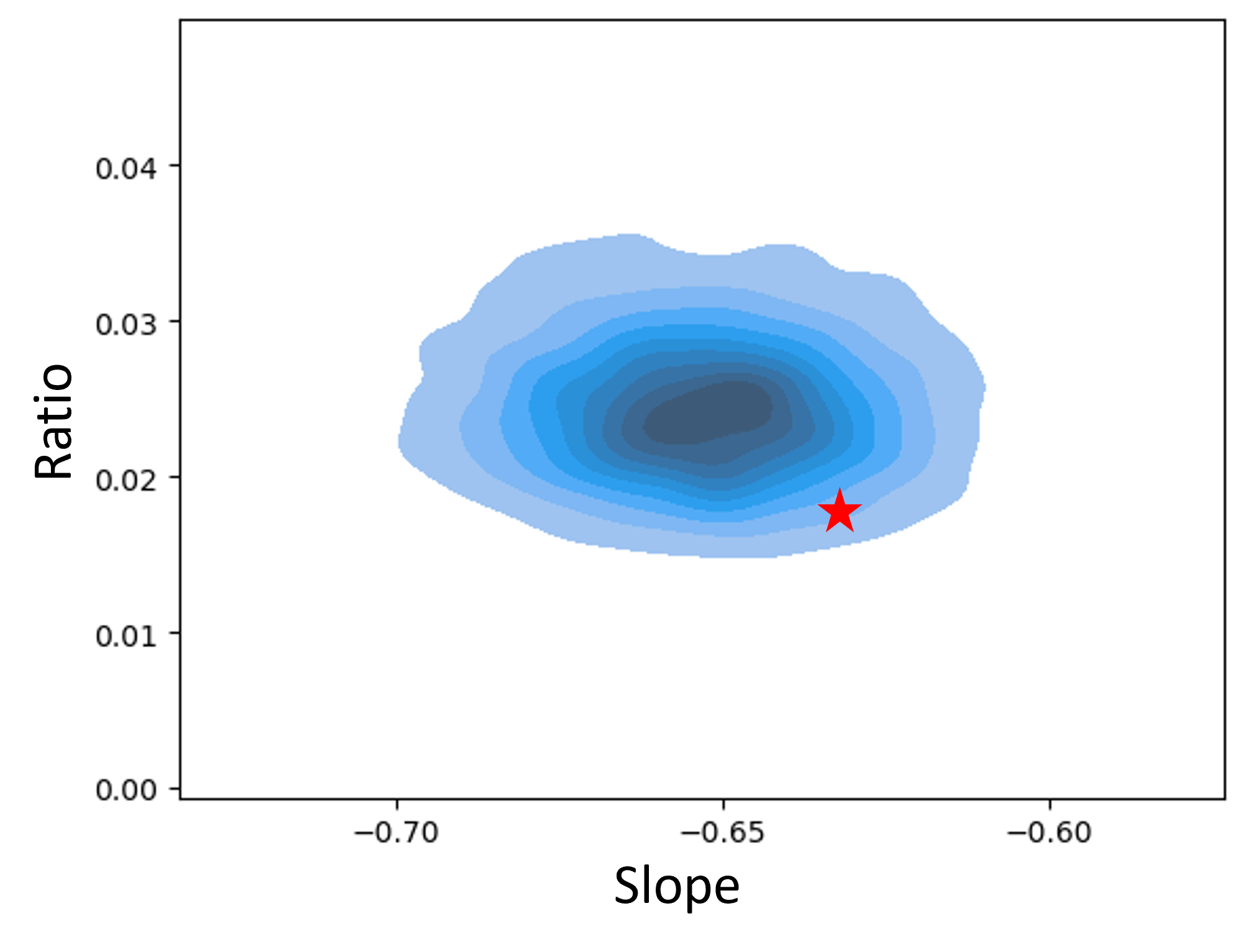}
\caption{}
\end{subfigure}
\caption{In both panels, the red star marks the slope of $\lambda_2$'s eigenvector and the ratio $\frac{\lambda_1}{\lambda_2}$ for the original movie dataset. Panel~(a): The blue circles represent the (slope, ratio) pairs from a principal component analysis of 2000 randomized movie crews. Panel~(b): A kernel density estimate of the (slope, ratio) pairs, verifying that the data in Panel~(a) approximately follow a bivariate normal distribution.}
\label{fig:ex1}
\end{figure}

From the 1496 unique individuals in the data, we randomly draw 2000 sets of 180 teams, each team consisting of 10 people, and then run a PCA analysis on each ensemble of 180 crews. We compute the eigenvalues $\lambda_1 < \lambda_2$ and eigenvectors of each covariance matrix and then compute the slope of the eigenvector belonging to $\lambda_2$ and the ratio $\frac{\lambda_1}{\lambda_2}$ of the eigenvalues. Plotting the $(slope, ratio)$ pairs in Figure~\ref{fig:ex1}, we see that they approximately follow a bivariate normal distribution, and we can use the PCA data results to create a kernel density estimate. Since the kernel density estimate approximates a continuous density function, it does not give exact probabilities at a single point. Instead, it tells us the relative likelihood of landing in a small region around that point. Calculating the probability that a point would land in a small neighborhood (defined here as an area with one tenth of the PCA data range) of our original data point, we find that it is about 1.9\%. Thus the point stemming from the original movie dataset is not a typical value expected from the bivariate normal distribution, and we conclude that the observed ensemble of teams was likely not the result of random selection.

\paragraph{Correlation with performance measures}

Finally, we provide two exploratory analyses of the correlation between $(\mathcal{D}$, $\mathcal{S})$ and team performance outcomes. 

\begin{figure}
    \centering
    \includegraphics[width = 0.6\textwidth]{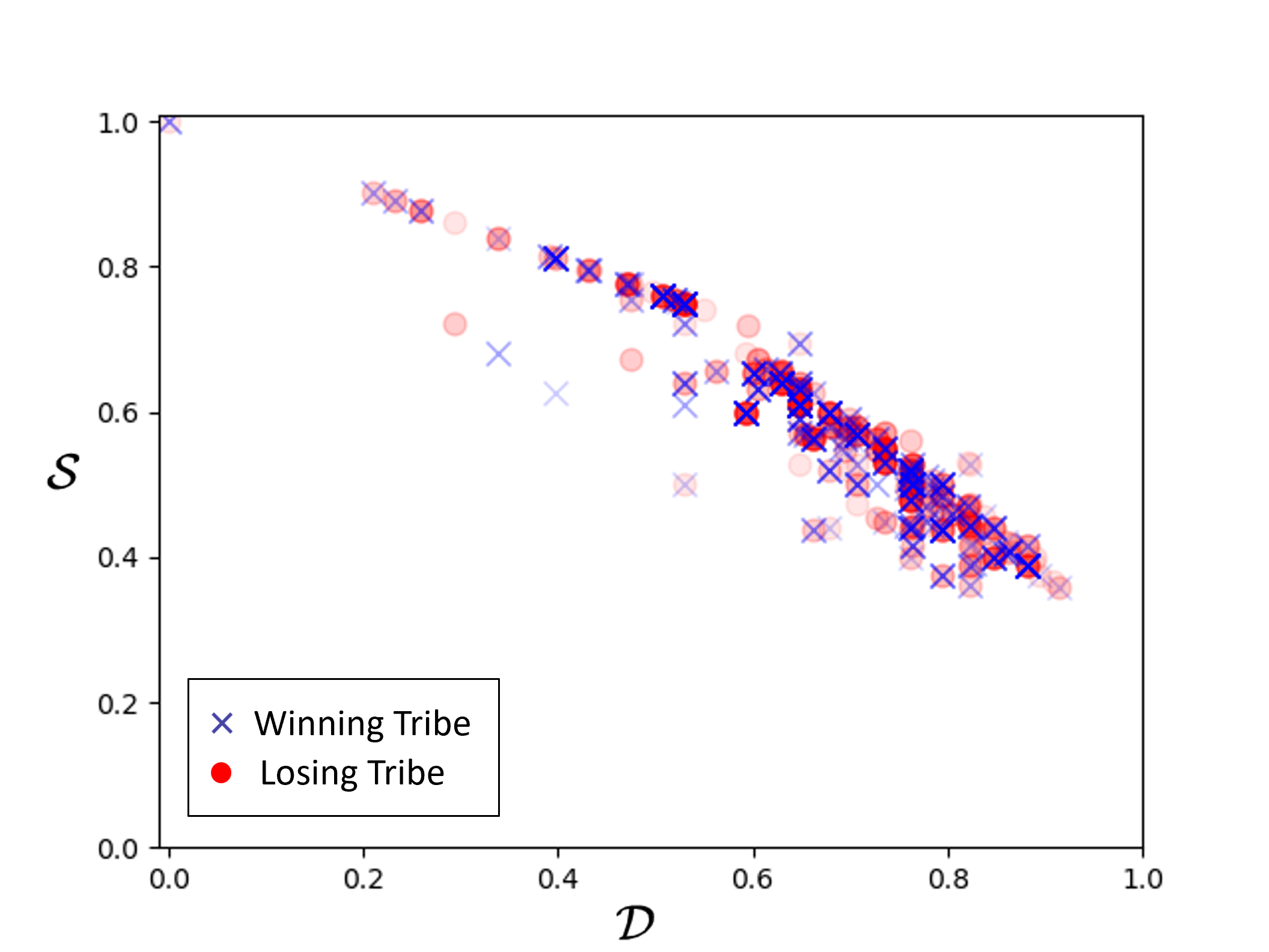}
    \caption{$(\mathcal{D},\mathcal{S})$ pairs for winning and losing teams in the \textit{Survivor} dataset considering gender and race/ethnicity ($T=2$). The blue X's represent winning outcomes and the large red circles represent losing outcomes. Darker shapes represent more teams.}
    \label{fig:S_outcomes}
\end{figure}

First, using the \textit{Survivor} data with the two traits gender and race/ethnicity, we explore whether higher team diversity led to more success in tribal challenges. Success of a team is defined as winning a challenge, while losing is coded as failure (in particular, we do not differentiate the quality of performance across winning teams). Figure~\ref{fig:S_outcomes} contains a scatter plot of the $(\mathcal{D},\mathcal{S})$ data in the $T=2$ case ($T=3$ results were similar), where winning and losing tribes are labeled differently. Visually, tribes with winning and losing outcomes have a similar spread in ($\mathcal{D},\mathcal{S}$) space. There is no evidence that the win/loss outcomes are related to ($\mathcal{D}$, $\mathcal{S}$) scores.

Second, for the North American company data, we use the real-valued industry-adjusted EBIT Margin as a measure of success (rather than a categorical win/loss measure), which allows us to compare how the relative magnitudes of $\mathcal{D}$ and $\mathcal{S}$ relate to relative magnitudes of performance. We employ an analysis from \cite{topazdiversity}, drawn from the theory that high Intersecting Diversity and high Shared Identity are both avenues to higher team performance. We test the hypothesis ``If B has better $\mathcal{D}$ \underline{and} $\mathcal{S}$ scores than A, then B has better performance than A." In their analysis, \cite{topazdiversity} examined demographics for U.S.\ states and used GDP per capita as the performance measure. They found this hypothesis to be true in 86 percent of the testable cases, along with evidence that this was not due to random chance. 

\begin{figure}
    \centering
    \begin{subfigure}[b]{0.49\textwidth}
        \centering
        \includegraphics[width=\textwidth]{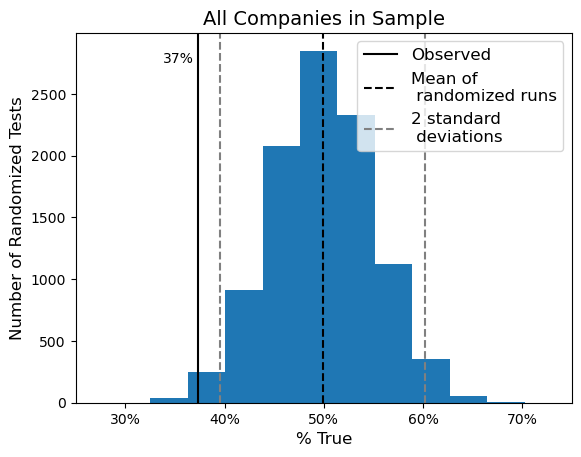} 
        \caption{}
        \label{fig:hypothesis_panel_a}
    \end{subfigure}
    \begin{subfigure}[b]{0.49\textwidth}
        \centering
        \includegraphics[width=\textwidth]{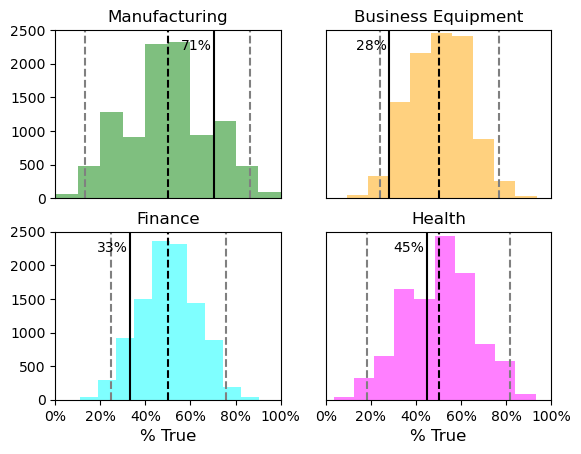}
        \caption{}
        \label{fig:hypothesis_panel_b}
    \end{subfigure}
    \caption{Percentage of true outcomes when testing the hypothesis ``Company B has better performance than Company A" in cases where Company B has both higher $\mathcal{D}$ and higher $\mathcal{S}$ scores than Company A. The histogram shows the results for the same companies with randomized rankings of performance, and the vertical line marks the result for the real performance data. Results in panel (a) include all companies and results in panel (b) are tests within industry. There were 1322 cases tested for all companies, 21 for Manufacturing, 31 for Health, 46 for Business Equipment, and 84 for Finance. Manufacturing's randomized results deviate slightly from a normal distribution shape due to a small sample size.}
    \label{fig:hypothesis_mainfig}
\end{figure}

We apply this analysis in our new context of North American companies. Here, for each two companies A and B for which B has both a higher Intersecting Diversity and higher Shared Identity score, we test the hypothesis ``B has a higher EBIT margin than A". We found 1322 cases where we were able to make a comparison, and the hypothesis was true in 37 percent of them, indicating that companies with higher scores in both $\mathcal{D}$ and $\mathcal{S}$ had lower EBIT margin scores on average. We checked whether this result was likely random by re-running the analysis 10,000 times using random shuffles of EBIT margin rankings to estimate a standard deviation. We found it to be statistically significant with a $p$-value choice of $0.05$ (see Figure \ref{fig:hypothesis_panel_a}), though running the same analysis using $\mathcal{S}_N$ instead of $\mathcal{S}$ did not find a statistically significant result. We also ran the analysis within the four largest industries in our sample (see Figure \ref{fig:hypothesis_panel_b}). Within these industries, there is no evidence that having higher $\mathcal{D}$ and $\mathcal{S}$ scores correlates with better financial performance in our sample. 

We stress that our analysis is exploratory and is not intended to study causal links between diversity, shared identity, and performance. It comes with many caveats. For the \textit{Survivor} data, we only inspect trends visually and do not provide a statistical test to determine significance or lack thereof. For the company data, we only examine correlations, and thus cannot determine whether company performance influences their leadership demographics, the reverse, or neither. The data are also likely not representative: we use a random sample of publicly-traded companies from the U.S.\ and Canada, exclude companies with incomplete data, and examine the unusual 2019-2023 period which may have been affected significantly by the COVID-19 pandemic. We cannot say conclusively that there is or is not a relationship between diversity, shared identity, and performance. However, we can conclude that we do not see the same evidence for the relationship between these metrics and high performance as was seen in \cite{topazdiversity} for a different application.


\section{Proofs of main results}
\label{sec:proofs}

In this section, we prove the results stated in  Section~\ref{sec:definitions}. 

\paragraph{Proof of Lemma~\ref{lemma:prob_int_S}.}

Let $G = \{ 1, 2, ..., N\}$ be a set of $N$ individuals equipped with the uniform probability measure $\mu$ defined by $\mu(U):=\frac{|U|}{N}$ for each $U\subset G$. Similarly, we equip the set $G\times G\setminus\Delta:=\{(i,j)\in G\times G\colon i\neq j\}$ with the measure $\nu$ defined by $\nu(U):=\frac{|U|}{N(N-1)}$ for $U\in G\times G\setminus\Delta$. We use
\[
\boldsymbol{\mathcal{C}} = \{ (c_1, c_2, ..., c_T)\colon 1 \leq c_t \leq v_t \text{ for } t \in 1, ..., T \} = \prod_{t=1}^T \{1,\ldots,v_t\}
\]
to denote the set of all possible aggregate identities, expressed as $T$-tuples where $c_t$ is the value for trait $t$. A function $Z\colon G\to\boldsymbol{\mathcal{C}},\; j\mapsto Z(j)=(c_1(j),\ldots,c_T(j))$ assigns to an individual $j$ their values for each trait $t$. Given such a function $Z$, we define the pushforward measure $p:=Z_\#\mu$ on $\boldsymbol{\mathcal{C}}$ so that
\[
p_c := p(c) = \mu(Z^{-1}(c)) = \frac{1}{N} \sum_{j\in G} \mathds{1}_{Z(j)=c} = P(Z=c)
\]
is the distribution of $Z$, where $\mathds{1}_\mathrm{statement}$ is $1$ if the statement is true and $0$ otherwise. Additionally, let $Y\colon G\times G\setminus\Delta\to\mathbb{N}$ be the function that assigns the number $Y(i,j)$ of shared traits to a pair $(i,j)$ of distinct individuals ($i\neq j$). Note that $Y = X|_O$ from the statement of Lemma~\ref{lemma:prob_int_S}, and we have
\begin{align}\label{e:Y}
Y\colon G\times G\setminus\Delta\longrightarrow \mathbb{N},\; (i,j)\longmapsto Y(i,j), \nonumber \\ Y(i,j) = \sum_{t=1}^T \mathds{1}_{Z_t(i)=Z_t(j)} =
\sum_{t=1}^T \sum_{v=1}^{v_t} \mathds{1}_{Z_t(i)=v} \mathds{1}_{Z_t(j)=v}.
\end{align}
The associated distribution $Y_\#\nu$ on $\mathbb{N}$ has the property that $P(Y=y):=Y_\#\nu(y)$ is the probability that two randomly selected individuals  sampled without replacement share exactly $y\in\mathbb{N}$ traits. Note that $0\leq Y(i,j)\leq T$ for all $(i,j)$ with $i \neq j$ and the support of $Y_\#\nu$ lies therefore in the set $\{0,\ldots,T\}\subset\mathbb{N}$. The expectation of $Y$ is given by
\begin{align}\label{e:EY}
E(Y) &= \sum_{y=0}^T y P(Y=y) = \int_\mathbb{N} y\,\mathrm{d}Y_\#\nu(y) \nonumber \\ &= 
\int_{G\times G\setminus\Delta} Y(i,j)\, \mathrm{d}\nu(i,j) = 
\frac{1}{N(N-1)} \sum_{(i,j)\in G\times G\setminus\Delta} Y(i,j).
\end{align}
Notice
\[
\mathcal{S}_N \stackrel{\eqref{eqn: OG S}}{=} 
\frac{\sum_{i>j} Y(i,j)}{T \binom{N}{2}}
= \frac{\sum_{i\neq j} Y(i,j)}{TN(N-1)}
= \frac{1}{TN(N-1)} \sum_{(i,j)\in G\times G\setminus\Delta} Y(i,j)
\stackrel{\eqref{e:EY}}{=} \frac{E(Y)}{T},
\]
which completes the proof of Lemma~\ref{lemma:prob_int_S}.

\paragraph{Proof of Theorem~\ref{main_thm}.}

We now turn to Theorem~\ref{main_thm} and first establish a relationship between $\mathcal{S}$ and $\mathcal{S}_N$. To do so, we consider \eqref{e:EY} and refine the expression by letting $i$ and $j$ become independent of each other:
\begin{align}\label{e:EYN}
E(Y) & = \frac{1}{N(N-1)} \sum_{i\neq j} Y(i,j) \stackrel{\eqref{e:Y}}{=}
\frac{1}{N(N-1)} \sum_{i\neq j} \sum_{t=1}^T \sum_{v=1}^{v_t} \mathds{1}_{Z_t(i)=v} \mathds{1}_{Z_t(j)=v} \\ \nonumber & =
\frac{1}{N(N-1)} \sum_{t=1}^T \sum_{v=1}^{v_t} \left( \sum_{i, j} \mathds{1}_{Z_t(i)=v} \mathds{1}_{Z_t(j)=v} - \sum_{j} \left( \mathds{1}_{Z_t(j)=v} \right)^2  \right) \\ 
\nonumber & = \frac{1}{N(N-1)} \sum_{t=1}^T \sum_{v=1}^{v_t} \left( \sum_{i} \left( \mathds{1}_{Z_t(i)=v} \left( \sum_{j} \mathds{1}_{Z_t(j)=v} \right) \right) - \sum_{j} (\mathds{1}_{Z_t(j)=v})^2 \right)
\\ \nonumber & = \frac{1}{N(N-1)} \sum_{t=1}^T \sum_{v=1}^{v_t} \left( \left( \sum_{i} \mathds{1}_{Z_t(i)=v} \right) \left( \sum_{j} \mathds{1}_{Z_t(j)=v} \right) - \sum_{j} \mathds{1}_{Z_t(j)=v} \right)
\\ \nonumber & =
\frac{1}{N(N-1)} \sum_{t=1}^T \sum_{v=1}^{v_t} \left( \left(\sum_{j} \mathds{1}_{Z_t(j)=v}\right)^2 - \sum_{j} \mathds{1}_{Z_t(j)=v} \right) \\ \nonumber & =
\frac{N}{N-1} \sum_{t=1}^T \sum_{v=1}^{v_t} \left(\sum_{j} \frac{1}{N} \mathds{1}_{Z_t(j)=v}\right)^2 - \frac{T}{N-1} \\ \nonumber & =
\frac{N}{N-1} \sum_{t=1}^T \sum_{v=1}^{v_t} P(Z_t=v)^2 - \frac{T}{N-1} =
\frac{N}{N-1} \sum_{t=1}^T \sum_{v=1}^{v_t} \left( \sum_{\underset{\scriptstyle c_t=v}{c\in\boldsymbol{\mathcal{C}}}} p_c\right)^2 - \frac{T}{N-1}.
\end{align}
Since
\begin{equation}\label{e:Sinfty}
\mathcal{S} \stackrel{\eqref{Sinf_def}}{:=} \frac{1}{T} \sum_{t=1}^T \sum_{v=1}^{v_t} \left( \sum_{\underset{\scriptstyle c_t=v}{c\in\boldsymbol{\mathcal{C}}}} p_c\right)^2,
\end{equation}
we see that
\[
\mathcal{S}_N = \frac{E(Y)}{T} = \frac{N}{N-1}\mathcal{S}-\frac{1}{N-1}
\]
or, alternatively, $\mathcal{S} = \left( 1 - \frac{1}{N}\right) \mathcal{S}_N +\frac{1}{N}$ as claimed.

Next, we can establish the bounds on $\mathcal{S}_N$ and $\mathcal{S}$ stated in Theorem~\ref{main_thm}(i). First, we calculate $P(Y=T)$. Since $Y=T$ holds precisely when the two sampled individuals share the same aggregate identity $c\in\boldsymbol{\mathcal{C}}$, we have $P(Y=T)=E(\tilde{Y})$ where $\tilde{Y}$ is the random variable corresponding to a single trait that takes on all possible aggregate identities $c\in\boldsymbol{\mathcal{C}}$ as values. Thus, proceeding as in \eqref{e:EYN} to calculate $E(\tilde{Y})$ for this single trait with $T=1$ and $c\in\boldsymbol{\mathcal{C}}$ in place of $v$, we see upon using the definition $\sum_{c\in\boldsymbol{\mathcal{C}}}p_c^2 \stackrel{\eqref{ID_defn}}{=}1-\frac{C-1}{C}\mathcal{D}$ that
\begin{equation}\label{e:PYT}
P(Y=T) = \frac{N}{N-1} \sum_{c\in\boldsymbol{\mathcal{C}}} p_c^2 - \frac{1}{N-1} =
\frac{N}{N-1} \left( 1-\frac{C-1}{C}\mathcal{D} \right) - \frac{1}{N-1}.
\end{equation}
Hence, we have
\[
\mathcal{S}_N = \frac{1}{T} \sum_{t=0}^T x P(Y=x) \geq P(Y=T) \stackrel{\eqref{e:PYT}}{=}
\frac{N}{N-1} \left( 1-\frac{C-1}{C}\mathcal{D} \right) - \frac{1}{N-1},
\]
and therefore obtain the lower bound
\[
\mathcal{S} \geq 1 - \frac{C-1}{C} \mathcal{D}.
\]
Furthermore, we see that
\begin{align*}
\mathcal{S}_N & =
\frac{1}{T} \sum_{x=0}^T x P(Y=x) = \frac{1}{T} \sum_{x=1}^T P(Y\geq x) =
\frac{1}{T} P(Y=T) + \frac{1}{T} \sum_{x=1}^{T-1} P(Y\geq x) \\ & \leq
\frac{1}{T} P(Y=T) + \frac{1}{T} \sum_{x=1}^{T-1} 1 \stackrel{\eqref{e:PYT}}{=}
\frac{1}{T} \left( \frac{N}{N-1} \left( 1-\frac{C-1}{C}\mathcal{D} \right) - \frac{1}{N-1} \right) + \frac{T-1}{T} \\ & \stackrel{\text{calculation}}{=}
\frac{N}{N-1} \left( 1-\frac{C-1}{TC}\mathcal{D} \right) - \frac{1}{N-1}
\end{align*}
and therefore obtain the upper bound
\[
\mathcal{S} \leq 1 - \frac{C-1}{TC} \mathcal{D}.
\]
This completes the proof of Theorem~\ref{main_thm}(i).

Next, we prove Theorem~\ref{main_thm}(ii). The expression
\[
\mathcal{S}(p) = \frac{1}{T} \sum_{t=1}^T \sum_{v=1}^{v_t} \left( \sum_{\underset{\scriptstyle c_t=v}{c\in\boldsymbol{\mathcal{C}}}} p_c\right)^2
\]
from \eqref{e:Sinfty} shows that $\mathcal{S}(p)$ is non-negative and quadratic in its argument $p:=(p_c)_{c\in\boldsymbol{\mathcal{C}}}$. We fix an index $\ell\in\boldsymbol{\mathcal{C}}$ and compute
\begin{align*}
\frac{\partial\mathcal{S}}{\partial p_\ell}(p)
& = \frac{2}{T} \sum_{t=1}^T \sum_{v=1}^{v_t} \left( \sum_{\underset{\scriptstyle c_t=v}{c\in\boldsymbol{\mathcal{C}}}} p_c \right) \left( \frac{\partial}{\partial p_\ell} \sum_{\underset{\scriptstyle c_t=v}{c\in\boldsymbol{\mathcal{C}}}} p_c \right)
= \frac{2}{T} \sum_{t=1}^T \sum_{v=1}^{v_t} \left( \sum_{\underset{\scriptstyle c_t=v}{c\in\boldsymbol{\mathcal{C}}}} p_c \right) \left( \sum_{\underset{\scriptstyle c_t=v}{c\in\boldsymbol{\mathcal{C}}}} \mathds{1}_{\ell=c} \right)
\\ & = \frac{2}{T} \sum_{t=1}^T \sum_{v=1}^{v_t} \left( \sum_{\underset{\scriptstyle c_t=v}{c\in\boldsymbol{\mathcal{C}}}} p_c \right) \mathds{1}_{\ell_t=v}
= \frac{2}{T} \sum_{t=1}^T \left( \sum_{\underset{\scriptstyle c_t=\ell_t}{c\in\boldsymbol{\mathcal{C}}}} p_c \right).
\end{align*}
We fix a second index $k\in\boldsymbol{\mathcal{C}}$ and calculate
\[
\frac{\partial^2\mathcal{S}}{\partial p_k \partial p_\ell}(p) =
\frac{2}{T} \sum_{t=1}^T \frac{\partial}{\partial p_k} \left( \sum_{\underset{\scriptstyle c_t=\ell_t}{c\in\boldsymbol{\mathcal{C}}}} p_c \right) =
\frac{2}{T} \sum_{t=1}^T \sum_{\underset{\scriptstyle c_t=\ell_t}{c\in\boldsymbol{\mathcal{C}}}} \mathds{1}_{c=k} =
\frac{2}{T} \sum_{t=1}^T \mathds{1}_{k_t=\ell_t}.
\]
Hence, $\mathcal{S}(p)=\frac12 p^* Q p$, where $Q$ is symmetric and positive semidefinite with entries $Q_{k\ell}=\frac{2}{T}\sum_{t=1}^T\mathds{1}_{k_t=\ell_t}\geq0$. Furthermore, the row sums of $Q$ are independent of $k$ and are given by
\[
\sum_{\ell\in\boldsymbol{\mathcal{C}}}Q_{k\ell} = 
\frac{2}{T} \sum_{t=1}^T \sum_{\ell\in\boldsymbol{\mathcal{C}}} \mathds{1}_{k_t=\ell_t} =
\frac{2}{T} \sum_{t=1}^T \prod_{\underset{\scriptstyle \tau\neq t}{1\leq\tau\leq T}} v_\tau =
\frac{2}{T} \sum_{t=1}^T \frac{C}{v_t} = 
\frac{2C}{T} \sum_{t=1}^T \frac{1}{v_t} =: \lambda_1.
\]
In particular, $\lambda_1$ is an eigenvalue of $Q$ with eigenvector $e=(1,\ldots,1)^* \in \mathbb{R}^C$.

We can now determine a second lower bound for $\mathcal{S}(p)$. Since $p$ corresponds to a probability measure, we have $\sum_{c\in\boldsymbol{\mathcal{C}}}p_c=1$ or, equivalently, $e^*p=1$. In particular, we have $(p-\frac{1}{C}e)\perp e$. Using that $Q$ is symmetric and positive semi-definite, we conclude that
\begin{eqnarray*}
\mathcal{S}(p) & = \frac{1}{2} p^* Q p =
\frac{\lambda_1}{2} \left|\frac{1}{C}e\right|^2 + \frac{1}{2} \underbrace{\left(p-\frac{1}{C}e\right)^t Q \left(p-\frac{1}{C}e\right)}_{\geq0} \geq
\frac{\lambda_1}{2} \left|\frac{1}{C}e\right|^2 =
\frac{\lambda_1}{2C} & \\ & =
\frac{1}{T} \sum_{t=1}^T \frac{1}{v_t} =: \mathcal{S}_\mathrm{min}.
\end{eqnarray*}
This completes the proof of Theorem~\ref{main_thm}(ii).

It remains to prove Theorem~\ref{main_thm}(iii). Recall that $\mathcal{D}(p)=\frac{C}{C-1}(1-|p|^2)$ and $\mathcal{S}(p)=\frac{1}{2}p^*Qp$ so that $\nabla\mathcal{D}(p)=-\frac{2C}{C-1}p$ and $\nabla\mathcal{S}(p)=Qp$. Hence,
\[
\langle\nabla\mathcal{D}(p),\nabla\mathcal{S}(p)\rangle =
-\frac{2C}{C-1} \langle p,Qp \rangle \leq -\frac{2C}{C-1} \lambda_1 |p|^2.
\]
Using $|p|^2=1-\frac{C-1}{C}\mathcal{D}\geq 1-\frac{C-1}{C}=\frac{1}{C}$, we obtain
\[
\langle\nabla\mathcal{D}(p),\nabla\mathcal{S}(p)\rangle \leq
-\frac{2\lambda_1}{C-1} = -\frac{4C}{C-1} \mathcal{S}_\mathrm{min} < 0
\]
as claimed.

\paragraph{Proof of Lemma~\ref{lemma:admissible}.}

Throughout the proof, we use the notation and results obtained above in this section, and we will also frequently refer to Figure~\ref{subfig:lemma2viz} for notation. Furthermore, for any given natural number $d$, we denote by $\{e_k\}_{1\leq k\leq d}$ the canonical basis vectors in $\mathbb{R}^d$.

We focus on exactly two traits ($T=2$) and assume, without loss of generality, that $1<v_1\leq v_2$. In particular, we have
\[
\boldsymbol{\mathcal{C}} = \{ (c_1, c_2)\colon 1 \leq c_t \leq v_t \text{ for } t \in 1,2 \} = \prod_{t=1}^2 \{1,\ldots,v_t\}.
\]
We denote by $\mathcal{P}$ the set of all probability measures $p=(p_c)_{c\in\boldsymbol{\mathcal{C}}}$ on $\boldsymbol{\mathcal{C}}$ so that $p_c\geq0$ and $\sum_{c\in\boldsymbol{\mathcal{C}}}p_c=1$ for each $p=(p_c)_{c\in\boldsymbol{\mathcal{C}}}\in\mathcal{P}$. We can then regard $\mathcal{D},\mathcal{S}:\mathcal{P}\to[0,1]$ as maps and are interested in describing the attainable set $\mathcal{A}:=\{(\mathcal{D},\mathcal{S})(p)\colon p\in\mathcal{P}\}$. We know that $P_1=(0,1)$ and $P_3=(1,\mathcal{S}_\mathrm{min})=(1,\frac{1}{2C}(v_1+v_2))$ are elements of $\mathcal{A}$, where $P_1$ is achieved by any Dirac measure, while $P_3$ is achieved by the uniform measure $p_c=\frac{1}{C}$ for all $c\in\boldsymbol{\mathcal{C}}$. 

We can represent $p\in\mathcal{P}$ by a matrix $(p_{ij})_{ij}\in\mathbb{R}^{v_1\times v_2}$, where $p_{ij}$ is the probability that the first trait has value $i$ and the second trait has value $j$. It will be more convenient to rearrange this matrix into a vector
\begin{equation}\label{e:x}
x = (r_1,\ldots,r_{v_1})^* \in \mathbb{R}^C = \mathbb{R}^{v_1v_2}, \qquad
r_i \in \mathbb{R}^{v_2} \mbox { for } i=1,\ldots,v_1,
\end{equation}
by concatenating the rows $r_i\in\mathbb{R}^{v_2}$ of $p$ for $i=1,\ldots,v_1$ (each row corresponds to a fixed value of the first trait) from top to bottom and then transpose the vector to form a column vector $x$. We then have $\mathcal{D}(x)=\frac{C}{C-1}(1-|x|^2)$ and $S(x)=\frac12 x^*Qx$. Using the representation \eqref{e:x} and the calculation of $Q$ established during the proof of Theorem~\ref{main_thm}, we see that $Q\in\mathbb{R}^{C\times C}$ is a $v_1\times v_1$ block matrix given by
\[
Q = \begin{pmatrix}
D       & 1_{v_2} & \ldots  & 1_{v_2} \\
1_{v_2} & D       & \ddots  & \vdots \\
\vdots  & \ddots  & \ddots  & 1_{v_2} \\
1_{v_2} & \ldots  & 1_{v_2} & D
\end{pmatrix} \in \mathbb{R}^{C\times C}, \qquad
D = \begin{pmatrix}
2 & 1 & \ldots & 1 \\
1 & 2 & \ddots & \vdots \\
\vdots & \ddots & \ddots & 1 \\
1 & \ldots & 1 & 2
\end{pmatrix} \in \mathbb{R}^{v_2\times v_2},
\]
where $1_{v_2}$ denotes the identity matrix in $\mathbb{R}^{v_2}$. The set $\mathcal{X}$ of vectors $x\in\mathbb{R}^C$ belonging to measures $p\in\mathcal{P}$ is given by $\mathcal{X}=\{x\in\mathbb{R}^C\colon e^*x=1,\; e_k^*x\geq0 \mbox{ for }k=1,\ldots,C\}$ where $e=(1, ..., 1)^*\in\mathbb{R}^C$.

We first calculate $L_{14}$. Since $|x|^2=1-\frac{C-1}{C}\mathcal{D}(x)\leq\mathcal{S}(x)$ by Theorem~\ref{main_thm}(i), it suffices to show that this lower bound for $\mathcal{S}$ can be achieved for fixed $|x|$. For $0\leq\alpha\leq\frac{v_1-1}{v_1}$, let $x(\alpha)=(1-\alpha)e_1+\frac{\alpha}{v_1-1}\sum_{i=1}^{v_1-1}e_{1+i(v_2+1)}$, then we can show that $x(\alpha)\in\mathcal{X}$ and $\mathcal{S}(x(\alpha))=|x(\alpha)|^2$ for $0\leq\alpha\leq\frac{v_1-1}{v_1}$. Furthermore, $x(0)=e_1$ and $x(\frac{v_1-1}{v_1})=\frac{1}{v_1}\sum_{i=0}^{v_1-1}e_{1+i(v_2+1)}$ with $P_1=(\mathcal{D},\mathcal{S})(x(0))=(1,0)$ and $P_4=(\mathcal{D},\mathcal{S})(x(\frac{v_1-1}{v_1}))=(\frac{C}{C-1}(1-\frac{1}{v_1}),\frac{1}{v_1})$ as claimed.

Before proceeding, we note that the symmetric matrix $Q$ has eigenvalues $\lambda=v_1+v_2, v_2, v_1, 0$ belonging to the pairwise orthogonal eigenspaces $X_m \subset \mathbb{R}^C$ for $m=v_1+v_2, v_2, v_1, 0$ of dimensions $1, v_1-1, v_2-1, (v_1-1)(v_2-1)$, respectively, which are given by
\begin{eqnarray*}
X_{v_1+v_2} & = & \mathrm{span}\{e\} \\
X_{v_2} & = & \mathrm{span}\{(\epsilon,-\epsilon,0,\ldots,0)^*, (0,\epsilon,-\epsilon,0,\ldots,0)^*, \dots, (0,\ldots,0,\epsilon,-\epsilon)^*\} \\
X_{v_1} & = & \mathrm{span}\{(z,\ldots,z)^*\colon \epsilon^*z=0, \; z\in\mathbb{R}^{v_2}\} \\
X_{0} & = &\mathrm{span}\{(z,-z,0,\ldots,0)^*, (0,z,-z,0,\ldots,0)^*, \dots, (0,\ldots,0,z,-z)^*\colon \epsilon^*z=0, \; z\in\mathbb{R}^{v_2}\},
\end{eqnarray*}
where $\epsilon=(1,\ldots,1)^*\in\mathbb{R}^{v_2}$. With $x^m\in X_m$ for $m=v_1+v_2,v_2,v_1,0$ and $x=\sum_m x^m\in\mathcal{X}$ so that $x^{v_1+v_2}=\frac{1}{C}e$, we have
\begin{eqnarray*}
       & \mathcal{D}(x) = \frac{C}{C-1} \left(1-\frac{1}{C}-|x^{v_2}|^2-|x^{v_1}|^2-|x^0|^2\right), & \\ &
\mathcal{S}(x) = \frac{1}{2} \left( \frac{v_1+v_2}{C} + v_2|x^{v_2}|^2 + v_1|x^{v_1}|^2 \right).
\end{eqnarray*}
In particular, we see that there is a point $Q_5=(q_5,\mathcal{S}_\mathrm{min})$ for some $q_5<1$ so that $\mathcal{S}$ does not change along the line segment $\overline{Q_5P_3}$: indeed, the vectors $x=\frac{1}{C}e+x^0$ with $x^0\in X_0$ lie in $\mathcal{X}$ as long as $e_k^*x^0\geq-\frac{1}{C}$ for $k=1,\ldots,C$, and we have $\mathcal{S}(x)=\mathcal{S}(\frac{1}{C}e)$ since $x^0$ lies in the null space of $Q$.

Next, we consider $L_{23}$. For a given value of $\mathcal{D}$ sufficiently close to one, we maximize $\mathcal{S}$ by evaluating at $x(\alpha)=\frac{1}{C}e+\alpha x^{v_2}$, where $x^{v_2}=((v_1-1)\epsilon,-\epsilon,\ldots,-\epsilon)^*\in X_{v_2}$ for $0\leq\alpha\leq\frac{1}{C}$. A calculation shows that
\[
L_{23} = \left\{ P_3 + \alpha^2\left(-\frac{C^2(v_1-1)}{C-1},\frac12 Cv_2(v_1-1) \right) \colon 0\leq\alpha\leq\frac{1}{C} \right\}
\]
and $L_{23}$ indeed ends at $P_2$ with $x(\frac{1}{C})=(\frac{1}{v_2}\epsilon,0,\ldots,0)^*$.

Finally, we claim that 
\[
L_{12} = \left\{ P_2 + \alpha^2 Cv_1(v_2-1) \left(-\frac{C}{C-1},\frac12 \right) \colon 0\leq\alpha\leq\frac{1}{C} \right\}.
\]
To see this, we note that this line segment connects $P_2$ for $\alpha=0$ with $P_1$ for $\alpha=\frac{1}{C}$ and has slope $\frac{-(C-1)}{2C}$. In particular, this line segment coincides with the upper bound for $\mathcal{S}$ obtained in Theorem~\ref{main_thm}(i), and it therefore suffices to show that each point on the segment is attained. The points in $L_{12}$ are attained at the elements $x(\alpha)$ for $0\leq\alpha\leq\frac{1}{C}$ given by
\begin{eqnarray*}
& x(\alpha) = \left( \frac{1}{v_2}\epsilon + \alpha v_1 z, 0,\ldots,0\right)^* = \frac{1}{C}(e+x^{v_2})
+ \alpha (x^{v_1} + x^{0}), & \\ &
z = (v_2-1,-1,\ldots,-1)^* \in\mathbb{R}^{v_2}, \quad
x^{v_1} = (z,\ldots,z)^*\in X_{v_1}, & \\ &
x^{0} = ((v_1-1)z,-z,\ldots,-z)^*\in X_0.
\end{eqnarray*}
This completes the proof of Lemma~\ref{lemma:admissible}.


\section{Discussion}
\label{sec:discussion}

The Intersecting Diversity ($\mathcal{D}$) and Shared Identity ($\mathcal{S}$) metrics provide unique quantitative tools to describe group demographics. 
Our work expanded on \cite{topazdiversity} to further examine the mathematical properties of these measures.  We provided probabilistic interpretations for these metrics: for two randomly chosen group members (sampled with replacement), $\mathcal{D}$ reflects the probability that they differ in at least one trait, while $\mathcal{S}$ is the expected percentage of shared traits. We defined $\mathcal{S}$ to be based on the original metric $\mathcal{S}_N$ from \cite{topazdiversity}, noting that $\mathcal{S}$ both aligns better with how $\mathcal{D}$ is defined and allows for direct comparisons between teams of varying sizes. The two measures $\mathcal{D}$ and $\mathcal{S}$ are anti-correlated: there is a trade-off between the two, and the set of possible values for $(\mathcal{D},\mathcal{S})$ pairs is restricted to an explicit polygonal shape in the unit square. To specify where ($\mathcal{D}, \mathcal{S}$) pairs will fall, we provided explicit bounds for their admissible region $\mathcal{R}$, which hold for any number of traits, and we further tightened these bounds by defining their attainable region $\mathcal{A}$ in the case where $T=2$. We also showed that $\mathcal{D}$ and $\mathcal{S}$ can be used to analyze datasets through hypothesis testing. In particular, we provided three illustrative case studies to quantify whether programs to increase diversity have been successful, test the null hypothesis that groups are drawn randomly from a demographically representative pool, and study the relationship between the composition of leadership teams and performance.

Our work has several limitations, many of which we reiterate from the discussion in \cite{topazdiversity}. First and foremost, any quantitative representation of human identity is inherently reductive, and metrics cannot capture the complexities of lived experience. Additionally, diversity is only part of the story regarding social justice and equality, and diversity is a social good regardless of its connection to team performance. Understanding this connection and the factors that influence effectiveness are tools to assist building diverse teams, not an argument against forming them. 

From an applied perspective, we assumed that all trait values are mutually exclusive, which is a simplification and often untrue for traits such as race/ethnicity in demographic datasets. One strength of the $\mathcal{D}$ and $\mathcal{S}$ measures is their ability to accommodate non-mutually exclusive trait values, so it would be very valuable to explore their properties with this assumption relaxed. It would also be useful to assess how our analysis would change if individuals are assigned a fractional value to each of their aggregate identities rather than forcing trait values to be binary at the individual level. This could preserve some of the measures' probabilistic properties, notably $\sum_{c \in \boldsymbol{\mathcal{C}}} p_c = 1$, but it may have other unforeseen effects. 

We also stressed the importance of context and care in choosing appropriate demographic values in one's analysis, but we did not systematically investigate how these choices influence the $\mathcal{D}$ and $\mathcal{S}$ metrics. Before using these measures in any sort of higher stakes context, it must be understood how trait composition, inclusion/omission, and structure affect the values of $\mathcal{D}$ and $\mathcal{S}$.
We showed in Section~\ref{sec:changes} that the individual properties of each trait influence how it contributes to the full $\mathcal{D}$ and $\mathcal{S}$ scores, but we did not conduct a full investigation. The sensitivity of $\mathcal{D}$ to trait value choices has been explored for a single axis of identity ($T = 1$) \cite{budescu2012measure}, but this analysis must be extended to account for multiple traits and must also be done for $\mathcal{S}$ and $(\mathcal{D}, \mathcal{S})$ pairs. Further work is needed to fully understand how the properties of each trait combine to  affect these multidimensional measures.


\section{Acknowledgements}
The authors gratefully acknowledge Chad Topaz and Heather Price for their thoughtful feedback on editing and paper's broader narrative. We extend special thanks to John Hand and his collaborators for generously sharing a random sample of their North American company data. Lastly, we thank Kevin Hu for insightful conversations and his contributions to an early proof of Lemma~1. This material is based upon work supported by the National Science Foundation (NSF). K. Slyman was supported by the NSF under grant DMS-2038039. B. Sandstede was partially supported by the NSF through grants DMS-2038039 and DMS-2106566. 


\section{Conflict of interest statement}
On behalf of all authors, the corresponding author states that there are no conflicts of interest.


\bibliographystyle{plain}
\bibliography{BIB}

\begin{thebibliography}{10}

\bibitem{bermiss2024does}
Sekou Bermiss, Jeremiah Green, and John~RM Hand.
\newblock Does greater diversity in executive race/ethnicity reliably predict better future firm financial performance?
\newblock {\em Journal of Economics, Race, and Policy}, 7(1):45--60, 2024.

\bibitem{boydstun2014importance}
Amber~E Boydstun, Shaun Bevan, and Herschel~F Thomas~III.
\newblock The importance of attention diversity and how to measure it.
\newblock {\em Policy Studies Journal}, 42(2):173--196, 2014.

\bibitem{budescu2012measure}
David~V Budescu and Mia Budescu.
\newblock How to measure diversity when you must.
\newblock {\em Psychological methods}, 17(2):215, 2012.

\bibitem{crenshaw2013demarginalizing}
Kimberl{\'e} Crenshaw.
\newblock Demarginalizing the intersection of race and sex: A black feminist critique of antidiscrimination doctrine, feminist theory and antiracist politics.
\newblock In {\em Feminist legal theories}, pages 23--51. Routledge, 2013.

\bibitem{mckinsey3-dolan2020diversity}
Kevin Dolan, Vivian Hunt, Sara Prince, and Sandra Sancier-Sultan.
\newblock Diversity still matters.
\newblock {\em McKinsey \& Company}, 2020.
\newblock \url{https://www.mckinsey.com/featured-insights/diversity-and-inclusion/diversity-still-matters}.

\bibitem{georgeac2023business}
Oriane~AM Georgeac and Aneeta Rattan.
\newblock The business case for diversity backfires: Detrimental effects of organizations' instrumental diversity rhetoric for underrepresented group members' sense of belonging.
\newblock {\em Journal of Personality and Social Psychology}, 124(1):69, 2023.

\bibitem{green2024mckinsey}
Jeremiah Green and John~RM Hand.
\newblock Mc{K}insey's diversity matters/delivers/wins results revisited.
\newblock {\em Econ Journal Watch}, 21(1), 2024.

\bibitem{herring2009does}
Cedric Herring.
\newblock Does diversity pay? {R}ace, gender, and the business case for diversity.
\newblock {\em American sociological review}, 74(2):208--224, 2009.

\bibitem{mckinsey1-hunt2015diversity}
Vivian Hunt, Dennis Layton, Sara Prince, et~al.
\newblock Diversity matters.
\newblock {\em McKinsey \& Company}, 1(1):15--29, 2015.
\newblock \url{https://www.mckinsey.com/capabilities/people-and-organizational-performance/our-insights/why-diversity-matters}.

\bibitem{mckinsey2-hunt2018delivering}
Vivian Hunt, Sara Prince, Sundiatu Dixon-Fyle, and Lareina Yee.
\newblock Delivering through diversity.
\newblock {\em McKinsey \& Company}, 2018.
\newblock \url{https://www.mckinsey.com/capabilities/people-and-organizational-performance/our-insights/delivering-through-diversity}.

\bibitem{jaccard1912distribution}
Paul Jaccard.
\newblock The distribution of the flora in the alpine zone.1.
\newblock {\em New phytologist}, 11(2):37--50, 1912.

\bibitem{joshi2009role}
Aparna Joshi and Hyuntak Roh.
\newblock The role of context in work team diversity research: A meta-analytic review.
\newblock {\em Academy of management journal}, 52(3):599--627, 2009.

\bibitem{kramer20011}
Roderick~M Kramer.
\newblock 1. organizational paranoia: Origins and dynamics.
\newblock {\em Research in organizational behavior}, 23:1--42, 2001.

\bibitem{milton2005identity}
Laurie~P Milton and James~D Westphal.
\newblock Identity confirmation networks and cooperation in work groups.
\newblock {\em Academy of Management Journal}, 48(2):191--212, 2005.

\bibitem{survivor_cran}
Daniel Oehm.
\newblock Survivor: Data from all seasons of survivor ({US}) {TV} series in tidy format.
\newblock {\em CRAN: Contributed Packages}, Apr 2021.
\newblock \url{http://dx.doi.org/10.32614/cran.package.survivoR}.

\bibitem{Variety}
Joe Otterson.
\newblock {CBS} sets diversity goal for unscripted shows - 50\% of cast to be people of color.
\newblock {\em Variety}, Nov 2020.
\newblock \url{https://variety.com/2020/tv/news/cbs-unscripted-shows-diversity-1234826214/}.

\bibitem{phillips2014diversity}
Katherine~W Phillips, Douglas Medin, Carol~D Lee, Megan Bang, Steven Bishop, and DN~Lee.
\newblock How diversity works.
\newblock {\em Scientific American}, 311(4):42--47, 2014.

\bibitem{roberto2015divergence}
Elizabeth Roberto.
\newblock The divergence index: A decomposable measure of segregation and inequality.
\newblock {\em arXiv preprint arXiv:1508.01167}, 2015.

\bibitem{ryabov_estimation_2022}
Alexey Ryabov, Bernd Blasius, Helmut Hillebrand, Irina Olenina, and Thilo Gross.
\newblock Estimation of functional diversity and species traits from ecological monitoring data.
\newblock {\em Proceedings of the National Academy of Sciences}, 119(43):e2118156119, October 2022.
\newblock Publisher: Proceedings of the National Academy of Sciences.

\bibitem{ryan2003stereotype}
Carey Ryan.
\newblock Stereotype accuracy.
\newblock {\em European review of social psychology}, 13(1):75--109, 2003.

\bibitem{sorenson1948method}
Th~Sorenson.
\newblock A method of establishing groups of equal amplitude in plant sociology based on similarity of species content, and its application to analysis of vegetation on danish commons.
\newblock {\em Kong Dan Vidensk Selsk Biol Skr}, 5:1--5, 1948.

\bibitem{stichman2010strength}
Amy~J Stichman, Kimberly~D Hassell, and Carol~A Archbold.
\newblock Strength in numbers? {A} test of {K}anter's theory of tokenism.
\newblock {\em Journal of Criminal Justice}, 38(4):633--639, 2010.

\bibitem{thomas1999cultural}
David~C Thomas.
\newblock Cultural diversity and work group effectiveness: An experimental study.
\newblock {\em Journal of cross-cultural psychology}, 30(2):242--263, 1999.

\bibitem{topazdiversity}
Chad~M Topaz, Heather~Z Brooks, Unchitta Kan, Bj{\"o}rn Sandstede, and Christian~Michael Smith.
\newblock Diversity, identity, and data.
\newblock {\em The American Mathematical Monthly}, 132(1):36--47, 2025.

\bibitem{verwijs2023double}
Christiaan Verwijs and Daniel Russo.
\newblock The double-edged sword of diversity: How diversity, conflict, and psychological safety impact software teams.
\newblock {\em IEEE Transactions on Software Engineering}, 50(1):141--157, 2023.

\bibitem{Deadline}
Peter White.
\newblock {CBS} sets diversity targets for reality casts; 50\% of talent must be {BIPOC} \& commits 25\% of unscripted development budget to {BIPOC} creatives.
\newblock {\em Deadline}, Nov 2020.
\newblock \url{https://deadline.com/2020/11/cbs-diversity-targets-reality-casts-bipoc-commits-unscripted-development-budget-1234611548/}.

\bibitem{wojcik_measuring_2025}
Laurie~Anne Wojcik, Ursula Gaedke, Ellen van Velzen, and Toni Klauschies.
\newblock Measuring overall functional diversity by aggregating its multiple facets: {Functional} richness, biomass evenness, trait evenness and dispersion.
\newblock {\em Methods in Ecology and Evolution}, 16(1):215--227, 2025.
\newblock \_eprint: https://onlinelibrary.wiley.com/doi/pdf/10.1111/2041-210X.14470.

\bibitem{zhu_trait_2017}
Linhai Zhu, Bojie Fu, Huoxing Zhu, Cong Wang, Lei Jiao, and Ji~Zhou.
\newblock Trait choice profoundly affected the ecological conclusions drawn from functional diversity measures.
\newblock {\em Scientific Reports}, 7(1):3643, June 2017.
\newblock Publisher: Nature Publishing Group.

\end{thebibliography}

\newpage

\end{document}